\documentclass[12pt,reqno]{amsart}

\usepackage{amsmath, amsfonts, amsthm, amssymb, color,  graphicx, mathrsfs, cite}
\usepackage{stmaryrd}

\makeatletter
\@addtoreset{figure}{section}
\makeatother

\theoremstyle{plain}
\newtheorem{Theorem}{Theorem}[section]
\newtheorem{Lemma}{Lemma}[section]
\newtheorem{Proposition}{Proposition}[section]

\theoremstyle{definition}

\theoremstyle{Remark}

\numberwithin{equation}{section}
\allowdisplaybreaks

\usepackage{ifpdf}
\ifpdf \usepackage[colorlinks=true, citecolor=blue, linkcolor=blue, urlcolor=blue]{hyperref} \fi

\newcommand{\RNum}[1]{\uppercase\expandafter{\romannumeral #1\relax}}

\newcommand{\ga}{\gamma}

\newcommand{\eps}{\epsilon}
\newcommand{\va}{\varepsilon}

\newcommand{\IT}{\int_0^t}

\begin{document}

\title[Rarefaction wave for the non-viscosity Navier-Stokes equations]
{Asymptotic stability of the rarefaction wave for
the non-viscous and heat-conductive ideal gas in half space}
\author{Meichen Hou$^{1,2}$}
\address{1.School of mathematical sciences, University of Chinese Academy of Sciences, Beijing,
100049, China}

\address{2. Institute of Applied Mathematics, AMSS, Beijing 100190, China Academy of Mathematics and Systems Science, Academia Sinica, Beijing 100190, China}

\email{meichenhou@amss.ac.cn}

\date{} \maketitle

\begin{abstract}This paper is concerned with the impermeable wall problem for an ideal polytropic model of non-viscous and heat-conductive gas in one-dimensional half space. It is shown that the 3-rarefaction wave
is stable under some smallness conditions. The proof is given by an elementary energy method and the key point is to control the boundary terms due to the less dissipativity of the system.

\bigbreak
 \noindent {\bf Keywords}: non-viscous; impermeable problem; rarefaction wave;
\end{abstract}

\section{Introduction}
In this paper we consider the one-dimensional initial boundary value problem (IBVP) for the equation of heat-conductive ideal gas without viscosity, which is modelled by following

\begin{equation}\label{1.1}
\left\{
\begin{aligned}
&\rho_{\tilde{t}}+(\rho u)_{\tilde{x}}=0,\\
&(\rho u)_{\tilde{t}}+(\rho u^2+p)_{\tilde{x}}=0,\\
&(\rho(e+\frac{u^2}{2}))_{\tilde{t}}+(\rho u(e+\frac{u^2}{2})+pu)_{\tilde{x}}=k\theta_{\tilde{x}\tilde{x}}.\\
\end{aligned}
\right.
\end{equation}
Here $(\tilde{t},\tilde{x})\in \mathbb{R}^+\times \mathbb{R}^+$ are Eulerian coordinates. And $\rho(\tilde{t},\tilde{x})>0,u(\tilde{t},\tilde{x}),e(\tilde{t},\tilde{x})>0,\theta(\tilde{t},\tilde{x})>0$ and $p(\tilde{t},\tilde{x})$ are density, fluid velocity, internal energy, absolute temperature and pressure respectively, while $k>0$ is the coefficient of the heat conduction, and it is a constant in this paper. We consider the ideal polytropic fluids so that

\begin{equation}\label{1.2}
p=R\rho\theta=A\rho^\ga\exp(\frac{\ga-1}{R}s),\quad e=\frac{R}{\ga-1}\theta,
\end{equation}
where $\ga>1, s$ are the adiabatic exponent of the gas and the specific entropy of the fluid respectively.
For the IBVP problem of (\ref{1.1}) in the half space $\tilde{x}>0$, the initial data is given by

\begin{equation}\label{1.3}
\left\{
\begin{array}{l}
(\rho_0,u_0,\theta_0)(\tilde{x})\triangleq (\rho,u,\theta)(0,\tilde{x})\rightarrow (\rho_+,u_+,\theta_+)
=z_+,\quad as \quad \tilde{x}\rightarrow +\infty,\\[2mm]
\inf_{\tilde{x}\in {\bf R_+}}(\rho_0,\theta_0)(\tilde{x})>0,
\end{array}
\right.
\end{equation}
where $\rho_+>0$, $\theta_+>0$ and $u_+\in
\mathbb{R}$ are given constants. Rewriting (\ref{1.1}) as following form

\begin{eqnarray}\label{1.4}
\left\{
\begin{array}{ll}
\rho_{\tilde{t}}+u\rho_{\tilde{x}}+\rho u_{\tilde{x}}=0,\\[2mm]
 u_t+uu_{\tilde{x}}+\frac{p}{\rho^2}\rho_{\tilde{x}}=-R\theta_{\tilde{x}},\\[2mm]
\frac{R}{\gamma-1}\theta_{\tilde{t}} -\frac{\kappa }{\rho}\theta_{\tilde{x}\tilde{x}}=-\frac{R}{\gamma-1}u \theta_{\tilde{x}}-R\theta u_{\tilde{x}},
\end{array}
\right.
\end{eqnarray}
one sees that the eigenvalues of  the hyperbolic part are
\begin{eqnarray}\label{1.5}
\widetilde{\lambda}_1=u-\widetilde{c}_s(\theta),\quad \widetilde{\lambda}_2=u+\widetilde{c}_s(\theta),
\end{eqnarray}
 where
 \begin{eqnarray}\label{1.6}
\widetilde{c}_s(\theta):=\sqrt{\frac{ p}{\rho}}=\sqrt{  R \theta}.
\end{eqnarray}

As pointed out by (\cite{MPPT}), the boundary conditions of (\ref{1.1}) depend on the sign of $\widetilde{\lambda}_1,\widetilde{\lambda}_2$. Because $\widetilde{\lambda}_1,\widetilde{\lambda}_2$ depend the solution itself, the boundary condition of (\ref{1.1}) should be proposed carefully. To simplify this,
the global solution $z=(\rho,u,\theta)$ of (\ref{1.1}) is considered in a small neighborhood $U(z_+)$ of the far field states $z_+=(\rho_+,u_+,\theta_+)$ such that the sign of $\widetilde{\lambda}_i(z)(i=1,2)$ is same as the sign of $\widetilde{\lambda}_i(z_+).$ Hence, the phase space is divided into following regions:
\begin{eqnarray*}
   & \widetilde{\Omega}^+_{sub}:=\left\{(\rho,u,\theta);\ 0<u<\widetilde{c}_s(\theta)\right\},\quad
     \widetilde{\Omega}^-_{sub}:=\left\{(\rho,u,\theta);\ -\widetilde{c}_s(\theta)<u<0\right\};\\[2mm]
     & \widetilde{\Omega}^+_{supper}:=\left\{(\rho,u,\theta);\ u>\widetilde{c}_s(\theta)\right\},\quad
     \widetilde{\Omega}^-_{supper}:=\left\{(\rho,u,\theta);\ u<-\widetilde{c}_s(\theta)\right\};\\[2mm]
    & \widetilde{\Gamma}^+_{trans}:=\left\{(\rho,u,\theta);\ u=\widetilde{c}_s(\theta)\right\},\quad
    \widetilde{\Gamma}^-_{trans}:=\left\{(\rho,u,\theta);\ u=-\widetilde{c}_s(\theta)\right\};\\[2mm]
    &\widetilde{\Gamma}^0_{sub}:=\left\{(\rho,u,\theta);\ u=0\right\}
   \end{eqnarray*}
and there are three cases for the boundary conditions of (\ref{1.1}):\\
Case (1):\ If  $z_+=(\rho_+,u_+,\theta_+) \in  \widetilde{\Omega}^-_{supper}$,  in the neighborhood of $U(z_+)$, $\widetilde{\lambda}_1(z)<0$,
$\widetilde{\lambda}_2(z)<0$, the boundary condition is
\begin{eqnarray}\label{1.7}
\theta(\tilde{t},0)=\theta_-.
\end{eqnarray}
Case (2):\ If  $z_+=(\rho_+,u_+,\theta_+) \in  \widetilde{\Omega}_{sub}^+\bigcup\widetilde{\Omega}_{sub}^-
\bigcup\widetilde{\Gamma}^0_{sub}$,  in the neighborhood of $U(z_+)$, $\widetilde{\lambda}_1(z)<0$,
$\widetilde{\lambda}_2(z)>0$, the boundary condition is
 \begin{eqnarray}\label{1.8}
u(\tilde{t},0)=u_-,\quad \theta(\tilde{t},0)=\theta_-.
\end{eqnarray}
Case (3):\ If  $z_+=(\rho_+,u_+,\theta_+) \in  \widetilde{\Omega}^+_{supper}$, in the neighborhood of $U(z_+)$, $\widetilde{\lambda}_1(z)>0$,
$\widetilde{\lambda}_2(z)>0$, the boundary condition is
 \begin{eqnarray}\label{1.9}
\rho(\tilde{t},0)=\rho_-,\quad u(\tilde{t},0)=u_-,\quad \theta(\tilde{t},0)=\theta_-.
\end{eqnarray}
Here $\rho_-,u_-,\theta_-$ are given constants. Above cases tell us that the IBVP problem of (\ref{1.1})
is very different from the IBVP problem of the viscous and heat-conductive ideal gas. For the latter, the system is 
\begin{equation}\label{1.10}
\left\{
\begin{aligned}
&\rho_{\tilde{t}}+(\rho u)_{\tilde{x}} = 0, \\
&(\rho u)_{\tilde{t}}+(\rho u^2+p)_{\tilde{x}} = \mu u_{\tilde{x}\tilde{x}}, \\
&(\rho (e+\frac{u^2}{2}))_{\tilde{t}}+(\rho u(e+\frac{u^2}{2})+pu)_{\tilde{x}} = k\theta_{\tilde{x}\tilde{x}}+(\mu uu_{\tilde{x}})_{\tilde{x}},
\end{aligned}
\right.
\end{equation}
where $\mu>0$ stands for the coefficient of viscosity. By \cite{Matsumura},  because the first equation is of hyperbolic, the second and third equations are of parabolic, then the boundary conditions of (\ref{1.10}) can be divided into three cases

\begin{equation}\label{1.11}
\left\{
\begin{aligned}
&\text{Case}\quad \RNum{1}: u(\tilde{t},0)=u_-<0,\quad\quad \theta(\tilde{t},0)=\theta_-,\\
&\text{Case}\quad \RNum{2}: u(\tilde{t},0)=u_-=0,\quad\quad \theta(\tilde{t},0)=\theta_-,\\
&\text{Case}\quad \RNum{3}: u(\tilde{t},0)=u_->0,\quad\quad \rho(\tilde{t},0)=\rho_-,\quad\quad \theta(\tilde{t},0)=\theta_-.
\end{aligned}
\right.
\end{equation}
Note that in Case $\RNum{1}-\RNum{2}$, the density $\rho$ on the boundary $\tilde{x}=0$ is unknown.
Usually, Case $\RNum{1}-\RNum{3}$ are called outflow problem ($u(\tilde{t},0)<0$), impermeable wall problem ($u(\tilde{t},0)=0$), inflow problem ($u(\tilde{t},0)>0$) respectively. Compare the boundary conditions ((\ref{1.7})-(\ref{1.9})) of (\ref{1.1}) with (\ref{1.11}), we immediately know that the IBVP problems of (\ref{1.1}) is more complicated than the related problems of (\ref{1.10}) due to the less dissipativity.

There are many works to the large time behavior of solutions to the Cauchy problem of systems (\ref{1.1}) and (\ref{1.10}), see (\cite{F-M},\cite{Huang-L-M}, \cite{Huang-Matsumura},\cite{Huang-Matsumura-Xin},
\cite{H-T-X},\cite{T.P.Liu-1},\cite{N-T-Z},\cite{Wang-Zhao}), all this works tell us that the large time behavior of solutions to the Cauchy problem is governed by the Riemann solutions to the corresponding compressible Euler equations.

For the IBVP problems of (\ref{1.1}) and (\ref{1.10}), the boundary layer maybe appear due to the boundary effect, see (\cite{Matsumura},\cite{M-Nishihara}). Therefore, it attracts many authors to study the related problems. For  the inflow problems of (\ref{1.10}), Qin-Wang (\cite{Qin-Wang},\cite{Qin-Wang-2011}) studied the superposition of multiple waves, i.e, rarefaction waves, boundary layer and contact wave. Recently, for the outflow problems of (\ref{1.1}), Nakamura-Nishibata \cite{Naka-Nishi1} studied the existence and stability of the boundary layer for the general symmetric hyperbolic-parabolic systems and their results contained the supersonic case of (\ref{1.1}), i.e., $z\in\widetilde{\Omega}_{supper}^-$ (Case (1)). And for the inflow problems of (\ref{1.1}), we have showned the existence of the boundary layer for $u_+>0$, moreover, the stability of the composition of rarefaction wave and boundary layer in supersonic case has been proved, i.e.,  $z\in\widetilde{\Omega}_{supper}^+$ (Case (3)), which will be appear, see (\cite{H-F}). 
For other related results, we refer to (\cite{F-H-W-Z},\cite{Huang-L-S}, \cite{H-M-S-1},\cite{H-M-S},\cite{Huang-Qin},\cite{F-Z}-\cite{Kawashima-Zhu},
\cite{M-Nishihara}-\cite{Nakamura2},\cite{Qin},\cite{Wan-Wang-Zou}) and some references therein.

For the impermeable wall problem, Matsumura-Mei \cite{Matsumura-Mei} studied the asymptotic stability of the viscous shock wave for the isentropic gas ($p=A\rho^\ga$). Matsumura-Nishihara {\cite{M-Nishihara1}} proved the stability of the rarefaction wave for the same model. Lately, Min-Qin \cite{Min-Qin} proved the asymptotic stability of the rarefaction wave of (\ref{1.10}) for some large perturbation. Motivated by this,  in this paper, we turn to study the stability of  the rarefaction wave for the impermeable wall problem of (\ref{1.1}). For this problem, $z_+\in\Omega_{sub}^+,z\in U(z_+),$
then $\widetilde{\lambda}_1(z)<0,\widetilde{\lambda}_2(z)>0$ immediately tell us that the boundary condition of this problem should be $(u,\theta)(t,0)=(u_-,\theta_-)$, here $u_-=0$(see Case (2)). Hence, the IBVP problem we consider here is (\ref{1.1}), (\ref{1.3}) and (\ref{1.8}), where $u_-=0$.
There are two main difficulties of this problem to overcome. The one is that both the density and the velocity equations are of hyperbolic, we need more higher order derivative estimates to recover the dissipativity of the hyperbolic part.
The other is that how to control the higher order derivatives of boundary terms (see ($H_2(\tau,0)$ in (\ref{3.51})). In order to solve the second difficulty, we use the relationship (\ref{3.56}) between the boundary terms in our estimates.

This paper is organized as follows. In Section 2, we do the transformation of coordinates for (\ref{1.1}) and  list some properties of the smooth approximation of rarefaction wave, then we state our main Theorem 2.1. In Section 3, we prove Theorem 2.1.
\vskip 0.5in
Notations. Throughout this paper, $c$ and $C$ denote some positive constants
(generally large).
For function spaces,~$L^p(\mathbb{R}_+)(1\leq p\leq \infty)$~denotes
the usual Lebesgue space on~$\mathbb{R}_+$~with norm~$\|{\cdot}\|_{L^p}$
and~$H^k(\mathbb{R}_+)$~the usual Sobolev space in the $L^2$ sense with norm~$\|\cdot\|_p$.
We note~$\|\cdot\|=\|\cdot\|_{L^2}$~for simplicity.
And $C^k(I; H^p)$ is the space of $k$-times continuously
differentiable functions on the interval $I$ with values in
$H^p(\mathbb{R}_+)$ and~$L^2(I; H^p)$ the space of~$L^2$-functions
on $I$ with values in~$H^p(\mathbb{R}_+)$.

\section{Preliminaries and Main results}
To simplify the system (\ref{1.1}), (\ref{1.3}) and (\ref{1.8}), we change the Eulerian coordinates $(\tilde{t},\tilde{x})$ into Lagrangian coordinates $(t,x)$
\begin{equation}\label{2.1}
t=\tilde{t},\quad\quad x=\int_{(0,0)}^{(\tilde{t},\tilde{x})}\rho(\tau,y)dy-(\rho u)(\tau,y)d\tau.
\end{equation}
Since $u_-=0,$ we let $v=\frac{1}{\rho}>0$ which is the specific volume of the gas, then the impermeable wall problem for (\ref{1.1}) becomes following
\begin{equation}\label{2.2}
\left\{
\begin{aligned}
&v_t-u_x=0, \quad\quad t>0,\quad x\in \mathbb{R}_+,\\
&u_t+p_x=0,\\
&(\frac{R}{\ga-1}\theta+\frac{u^2}{2})_t+(pu)_x=k(\frac{\theta_x}{v})_x,\\
&(v_0,u_0,\theta_0)(x)\triangleq (v,u,\theta)(0,x)\rightarrow (v_+,u_+,\theta_+)\quad\quad x\rightarrow +\infty, \\
&\inf_{x\in \mathbb{R}_+}(v_0,\theta_0)(x)>0,\quad u(t,0)=0,\quad\theta(t,0)=\theta_-,
\end{aligned}
\right.
\end{equation}
where 

\begin{equation}\label{2.2-1}
p=p(v,\theta)=\frac{R\theta}{v}
=A{v}^{-\ga}\exp(\frac{\ga-1}{R}s)
=p(v,s),
\end{equation}
and the initial data satisfy $u_0(0)=0,\theta_0(0)=\theta_-$ as compatibility conditions.
We have known that the corresponding hyperbolic system of (\ref{2.2}) has three characteristic speeds
\begin{equation}\label{2.3}
\lambda_1(v,\theta)=-\frac{\sqrt{R\gamma\theta}}{v}<0,\quad\quad \lambda_2=0,\quad\quad \lambda_3(v,\theta)=\frac{\sqrt{R\gamma\theta}}{v}>0.
\end{equation}

Now we turn to list some properties of the 3-rarefaction wave. The 3-rarefaction wave curve through the right-hand side state $(v_+,u_+,\theta_+)$
is
\begin{eqnarray}\label{2.4-1}
R_3(v_+,u_+,\theta_+):&=&\big\{(v,u,\theta): v>v_+, v^{\gamma-1}\theta=v^{\gamma-1}_+\theta_+,\\[2mm]
&&u=u_+-\int^v_{v_+}\sqrt{R\gamma v_+^{\gamma-1}\theta_+}\xi^{-\frac{\gamma+1}{2}}d\xi\big\}.
\end{eqnarray}
Precisely to say, if $v_+,u_+,\theta_+,\theta_-$ satisfy following condition
\begin{eqnarray}\label{2.4}
0<u_+=\int^{(\frac{\theta_+}{\theta_-})^{\frac{1}{\gamma-1}}v_+}_{v_+}\sqrt{R\gamma v_+^{\gamma-1}\theta_+}\xi^{-\frac{\gamma+1}{2}}d\xi,
\end{eqnarray}
then the 3-rarefaction wave $(v^r,u^r,\theta^r)(\frac{x}{t})$  is the unique weak solution which is global in time
to the following Riemann problem

\begin{equation}\label{2.5}
\left\{
\begin{aligned}
&v_t^r-u_x^r=0, \\
&u_t^r+p_x^r=0,  \\
&(\frac{R}{\ga-1}\theta^r+\frac{u^{r^2}}{2})_t+(p^ru^r)_x=0, \\
&(v^r,u^r,\theta^r)(t,0)=
\begin{cases}
&(v_-,0,\theta_-),\quad x<0, \\
&(v_+,u_+,\theta_+),\quad x>0.
\end{cases}
\end{aligned}
\right.
\end{equation}
Here $p^r=\frac{R\theta^r}{v^r}$, $v_->v_+,\theta_-<\theta_+,$ and $v_-=(\frac{\theta_+}{\theta_-})^{\frac{1}{\ga-1}}v_+.$ To study the large time behavior
of  the solutions to the impermeable wall problem (\ref{2.2}), we
construct a smooth approximation $(\tilde{v},\tilde{u},\tilde{\theta})(t,x)$ of 3-rarefaction wave $(v^r,u^r,\theta^r)(\frac{x}{t})$. Same as \cite{H-M-S-2}, we consider the following Cauchy problem
\begin{eqnarray}\label{2.6}
\left\{
\begin{array}{ll}
w_t+ww_x=0,\quad\quad\\[2mm]
w(x,0)=\left\{
\begin{array}{ll}
w_-, \quad x<0,\\[2mm]
w_-+C_q\delta^r \int^{\epsilon x}_0y^qe^{-y}dy, \quad x>0.
\end{array}
\right.
\end{array}
\right.
\end{eqnarray}
Here $\delta^r=w_+-w_->0$, $0<\epsilon\leq 1$ is a constant which will be determined later and $q\geq 10, C_q$ are two constants such that $C_q \int^{\infty}_0y^qe^{-y}dy=1$.
Then the smooth approximation $(\tilde{v},\tilde{u},\tilde{\theta})(t,x)$ is constructed in the following way
\begin{eqnarray}\label{2.7}
\left\{
\begin{array}{ll}
(\frac{\sqrt{R\gamma \tilde{\theta}}}{\tilde{v}})(t,x)=w(1+t,x),\\[2mm]
(\tilde{v}^{\gamma-1}\tilde{\theta})(t,x)=v^{\gamma-1}_+\theta_+,\quad x \in \mathbb{R},\quad t>0.\\[2mm]
\tilde{u}=u_+-\int^{\tilde{v}}_{v_+}\sqrt{R\gamma v^{\gamma-1}_+\theta_+}\xi^{-\frac{\gamma+1}{2}}d\xi.
\end{array}
\right.
\end{eqnarray}
Because $\lambda_3(v_-,\theta_-)>0$, we immediately know that both the 3-rarefaction wave $(v^r,u^r,\theta^r)$ and
its smooth approximation $(\tilde{v},\tilde{u},\tilde{\theta})$ are constants on $(t,x)\in(0,+\infty)\times R_-.$
Then we use $(\tilde{v},\tilde{u},\tilde{\theta})(t,x)$  to represent $(\tilde{v},\tilde{u},\tilde{\theta})(t,x) |_{x\geq0}$, we have

\begin{equation}\label{2.8}
\left\{
\begin{aligned}
&\tilde{v}_t-\tilde{u}_x=0, \\
&\tilde{u}_t+\tilde{p}_x=0,\quad\quad t>0,\quad x\in \mathbb{R}_+  \\
&(\frac{R}{\ga-1}\tilde{\theta}+\frac{\tilde{u}^{2}}{2})_t+(\tilde{p}\tilde{u})_x=0, \\
&(\tilde{v},\tilde{u},\tilde{\theta})(t,0)=(v_-,0,\theta_-),\\
&(\tilde{v}_0,\tilde{u}_0,\tilde{\theta}_0)(x)\rightarrow
\begin{cases}
&(v_-,0,\theta_-),\quad\quad x\rightarrow 0_+, \\
&(v_+,u_+,\theta_+),\quad\quad x\rightarrow +\infty,
\end{cases}
\end{aligned}
\right.
\end{equation}
where $\tilde{p}=p(\tilde{v},\tilde{\theta})=\frac{R\tilde{\theta}}{\tilde{v}}$.
We get following Lemma.

\begin{Lemma}\label{L1}(Smooth rarefaction wave)(\cite{H-M-S-2})$(\tilde{v},\tilde{u},\tilde{\theta})(t,x)$ satisfies

(1)$0\leq -\tilde{v}_x(t,x),\frac{R}{\ga-1}\tilde{\theta}_x(t,x)\leq C\tilde{u}_x(t,x)$;

(2)For any p($1\leq p\leq +\infty$), there exists a constant $C_{pq}$ such that
\begin{equation}\label{2.17}
\begin{aligned}
&\|(\tilde{v}_x,\tilde{u}_x,\tilde{\theta}_x)(t)\|_{L^p}\leq C_{pq}\min\{\delta^r\eps^{1-\frac{1}{p}},
(\delta^r)^\frac{1}{p}(1+t)^{-1+\frac{1}{p}}\}, \\
&\|(\tilde{v}_{xx},\tilde{u}_{xx},\tilde{\theta}_{xx})(t)\|_{L^p} \leq C_{pq}\min
\{\delta^r\eps^{2-\frac{1}{p}},((\delta^r)^{\frac{1}{p}}+(\delta^r)^{\frac{1}{q}})(1+t)^{-1+\frac{1}{q}}\},
\end{aligned}
\end{equation}

(3)If $x<\lambda_3(v_-,\theta_-)(1+t),$ then
$(\tilde{v},\tilde{u},\tilde{\theta})(t,x)\equiv(v_-,u_-,\theta_-),$

(4)$\lim_{t\rightarrow+\infty}\sup_{\xi\in R_+}|(\tilde{v},\tilde{u},\tilde{\theta})(t,x)-(v^r,u^r,\theta^r)(\frac{x}{1+t})|=0.$
\end{Lemma}

We assume that the initial data satisfy
\begin{equation}\label{2.18}
(v_0-\tilde{v}_0,u_0-\tilde{u}_0,\theta_0-\tilde{\theta}_0)\in H^2(\mathbb{R}_+),\quad
(v_t-\tilde{v}_t,u_t-\tilde{u}_t,\theta_t-\tilde{\theta}_t)(0,x)\in H^1(\mathbb{R}_+).
\end{equation}
Then the stability of 3-rarefaction wave is listed as follows

\begin{Theorem}\label{t1} Assume that $\ga>1$ and
the relationship between  $v_+,u_+,\theta_+,$ and $\theta_-$ satisfy
\begin{equation}\label{2.19}
u_+=\int^{(\frac{\theta_+}{\theta_-})^{\frac{1}{\gamma-1}}v_+}_{v_+}\sqrt{R\gamma v_+^{\gamma-1}\theta_+}\xi^{-\frac{\gamma+1}{2}}d\xi,\quad\quad \theta_+-\theta_->0.
\end{equation}
If the initial data satisfy (\ref{2.18}),
then there exist positive constants $\eps_1,\va_0$ such that if  $\eps\leq\eps_1$ and
\begin{equation}\label{2.20}
\|v_0-\tilde{v}_0,u_0-\tilde{u}_0,\theta_0-\tilde{\theta}_0\|_2+\|(v_t-\tilde{v}_t,u_t-\tilde{u}_t,\theta_t-\tilde{\theta}_t)(0)\|_1\leq \va_0,
\end{equation}
 the impermeable wall problem (\ref{2.2}) has a unique global solution $(v,u,\theta)(t,x)$ which satisfy
\begin{equation}\label{2.21}
\left\{
\begin{aligned}
&(v-\tilde{v},u-\tilde{u},\theta-\tilde{\theta})\in C([0,+\infty);H^2(R_+)),
(v-\tilde{v},u-\tilde{u},\theta-\tilde{\theta})_t\in C([0,+\infty);H^1(R_+)),  \\
&(v-\tilde{v},u-\tilde{u})_x\in L^2(0,+\infty;H^1(R_+)),(\theta-\tilde{\theta})_x\in L^2(0,+\infty;H^2(R_+)), \\
&(\theta-\tilde{\theta})_t\in L^2(0,+\infty;H^2(R_+)).
\end{aligned}
\right.
\end{equation}
Moreover,
\begin{equation}\label{2.27}
\sup_{x\geq 0}|(v,u,\theta)(t,x)-(v^r,u^r,\theta^r)(\frac{x}{t})|\rightarrow 0, \quad as
\quad t\rightarrow +\infty.
\end{equation}
\end{Theorem}

\section{Proof of Theorem 2.1}
This section we mainly proof Theorem \ref{t1}. In Section 3.1, we shall reformulate system (\ref{2.2}) to a new perturbed system and show the local existence of the solution, see  Proposition 3.1. In Section 3.2, a priori estimates will be investigated, see Proposition 3.2. Finally ,we will combine this two Propositions to
get the stability of the 3-rarefaction wave.

\subsection{Local existence of the solution}
Now we turn to reformulate the system (\ref{2.2}). Put the perturbation $(\phi,\psi,\xi)(t,x)$ by
\begin{equation}\label{3.1}
(\phi,\psi,\xi)(t,x)=(v,u,\theta)(t,x)-(\tilde{v},\tilde{u},\tilde{\theta})(t,x),
\end{equation}
then the reformulated problem is
\begin{equation}\label{3.2}
\left\{
\begin{aligned}
&\phi_t-\psi_x=0,\quad\quad t>0,\quad x\in \mathbb{R}_+\\
&\psi_t+(\frac{R\xi-\tilde{p}\phi}{v})_x=0,\\
&\frac{R}{\gamma-1}\xi_t+p\psi_x+\tilde{u}_x(p-\tilde{p})
=k(\frac{\theta_x}{v})_x,\\
&(\psi,\xi)(t,0)=(0,0), \\
&(\phi,\psi,\xi)(0,x)\triangleq(\phi_0,\psi_0,\xi_0)(x)\rightarrow (0,0,0),\quad as \quad x\rightarrow +\infty.
\end{aligned}
\right.
\end{equation}
Define the solution space as
\begin{equation}\label{3.3}
\begin{aligned}
&\mathbb{X}_{m_1,m_2,M}(0,T):=\bigg\{(\phi,\psi,\xi)\in C([0,T];H^2(\mathbb{R}_+)),
(\phi,\psi,\xi)_t\in C([0,T];H^1(\mathbb{R}_+)),\\
&(\phi,\psi)_x\in L^2(0,T;H^1(\mathbb{R}_+)),\xi_x\in L^2(0,T;H^2(\mathbb{R}_+)),
\xi_t\in L^2(0,T;H^2(\mathbb{R}_+)),\\
&\text{with}\quad\sup_{[0,T]}\|(\phi,\psi,\xi)(t)\|_2+\|(\phi_t,\psi_t,\xi_t)(t)\|_1\leq M,
\inf_{\mathbb{R}_+\times[0,T]}(\tilde{v}+\phi)(t,x)\geq m_1,\\
&\inf_{\mathbb{R}_+\times[0,T]}(\tilde{\theta}+\xi)(t,x)\geq m_2\bigg\}
\end{aligned}
\end{equation}
for some positive constants $m_1,m_2,M$.

\begin{Proposition}\label{p1} (Local existence) For any given initial data $(\phi_0,\psi_0,\xi_0)\in H^2(\mathbb{R}_+)$ and $(\phi_t,\psi_t,\xi_t)(0,x)\in H^1(\mathbb{R}_+)$,
there exists positive constant $M_0$ such that if
$\|(\phi_0,\psi_0,\xi_0)\|_2+\|(\phi_t,\psi_t,\xi_t)(0)\|_1\leq M(\tilde{C}M\leq M_0)$and $\inf_{\mathbb{R}_+\times[0,T]}(\tilde{v}+\phi_0)\geq \frac{3v_+}{8},\inf_{\mathbb{R}_+\times[0,T]}(\tilde{\theta}+\xi_0)\geq \frac{3\theta_-}{8}$, then there exists $t_0=t_0(M_0)>0$
such that (\ref{3.2}) has a unique solution $(\phi,\psi,\xi)\in\mathbb{X}_{\frac{3v_+}{8},\frac{3\theta_-}{8},\tilde{C}M}(0,t_0).$
\end{Proposition}

Proof. Rewrite $(\ref{3.2})$ with following forms,

\begin{equation}\label{3.4}
\left\{
\begin{aligned}
&\phi_t-\psi_x=0,\\
&\psi_t-\frac{p}{v}\phi_x
=g_1:=g_1(\phi,\xi,\xi_x),\\
&\frac{R}{\gamma-1}\xi_t-k\frac{\xi_{xx}}{\tilde{v}+\phi}=g_2:=g_2(\phi,\xi,\phi_x,\psi_x,\xi_x),\\
&\psi(t,0)=0,\xi(t,0)=0,\\
&(\phi(0,x),\psi(0,x),\xi(0,x))=(\phi_0(x),\psi_0(x),\xi_0(x))\rightarrow (0,0,0),\quad as \quad x\rightarrow +\infty,
\end{aligned}
\right.
\end{equation}
where

\begin{equation}\label{3.5}
\begin{aligned}
&g_1(\phi,\xi,\xi_x)=\frac{\tilde{p}_x\phi}{\tilde{v}+\phi}-\frac{\tilde{p}\tilde{v}_x\phi}
{(\tilde{v}+\phi)^2}-\frac{R\xi_x}{\tilde{v}+\phi}+\frac{R\xi\tilde{v}_x}{(\tilde{v}+\phi)^2},\\
&g_2(\phi,\xi,\phi_x,\psi_x,\xi_x)=-k\frac{\xi_x(\tilde{v}+\phi)_x}{(\tilde{v}+\phi)^2}
-p\psi_x-\tilde{u}_x(p-\tilde{p})+k(\frac{\tilde{\theta}_x}{\tilde{v}+\phi})_x.
\end{aligned}
\end{equation}

Now we approximate the initial data $(\phi_0,\psi_0,\xi_0)$ by 
$ (\phi_{0j},\psi_{0j},\xi_{0j})\in H^m\cap H^2$, $m\geq 6$ such that
\begin{equation}\label{3.6}
(\phi_{0j},\psi_{0j},\xi_{0j})\rightarrow (\phi_0,\psi_0,\xi_0)\quad \text{strongly}\quad\quad \text{in}\quad\quad H^m
\end{equation}
as $j\rightarrow\infty$ and 
\begin{equation}\label{3.6-1}
 \|(\phi_{0j},\psi_{0j})\|_2+\|\xi_{0j}\|_3\leq\frac{3}{2}M, 
\end{equation}
furthermore, $\inf_{\mathbb{R}_+}(\tilde{v}+\phi_{0j})(t,x)\geq \frac{3v_+}{8}, \inf_{\mathbb{R}_+}(\tilde{\theta}+\xi_{0j})(t,x)\geq \frac{3\theta_-}{8}$ hold for any $j\geq 1.$ (Here because our initial data satisfy $\|(\phi_0,\psi_0,\xi_0)\|_2+\|(\phi_t,\psi_t,\xi_t)(0)\|_1\leq M$, by system 
(\ref{3.2}), we know that the norm is equivalent to $\|(\phi_0,\psi_0)\|_2+\|\xi_0\|_3$, that's why we construct approximated sequence satisfy (\ref{3.6-1})).

We will use the iteration method to prove our Proposition 3.1.  Define the sequence
$\{(\phi^{(n)}_j(t,x),\psi^{(n)}_j(t,x),\xi^{(n)}_j(t,x))\}$ for each $j$ so that

\begin{equation}\label{3.7}
(\phi^{(0)}_j,\psi^{(0)}_j,\xi^{(0)}_j)(t,x)=(\phi_{0j},\psi_{0j},\xi_{0j})(x),
\end{equation}
and for a given $(\phi^{(n-1)}_j,\psi^{(n-1)}_j,\xi^{(n-1)}_j)(t,x),$
$(\phi^{(n)}_j,\psi^{(n)}_j,\xi^{(n)}_j)(t,x)$ is a solution to following equation

\begin{equation}\label{3.8}
\left\{
\begin{aligned}
&\phi_{jt}^{(n)}-\psi_{jx}^{(n)}=0\\
&\psi_{jt}^{(n)}-(\frac{R(\tilde{\theta}+\xi^{(n-1)}_j)}{(\tilde{v}+\phi_j^{(n-1)})^2})\phi^{(n)}_{jx}
=g_1^{(n-1)}=g_1^{(n-1)}(\phi^{(n-1)}_j,\xi^{(n-1)}_j,\xi^{(n-1)}_{jx})\\
&\frac{R}{\gamma-1}\xi_{jt}^{(n)}-k\frac{\xi_{jxx}^{(n)}}{\tilde{v}+\phi_{j}^{(n-1)}}=g_2^{(n-1)}
=g_2^{(n-1)}(\phi_{j}^{(n-1)},\xi_j^{(n-1)},\phi_{jx}^{(n-1)},\psi_{jx}^{(n-1)},\xi_{jx}^{(n-1)})\\
&\psi^{(n)}_j(t,0)=0,\quad\xi^{(n)}_j(t,0)=0\\
&(\phi^{(n)}_j,\psi^{(n)}_j,\xi^{(n)}_j)(0,x)=(\phi_{0j},\psi_{0j},\xi_{0j})(x),
\end{aligned}
\right.
\end{equation}

where
\begin{equation}\label{3.9}
\begin{aligned}
&g_1^{(n-1)}(\phi^{(n-1)}_j,\xi^{(n-1)}_j,\xi^{(n-1)}_{jx})\\
&=\frac{\tilde{p}_x\phi^{(n-1)}}{\tilde{v}+\phi^{(n-1)}_j}-\frac{\tilde{p}\tilde{v}_x\phi^{(n-1)}_j}
{(\tilde{v}+\phi^{(n-1)}_j)^2}-\frac{R\xi^{(n-1)}_{jx}}{\tilde{v}+\phi^{(n-1)}_j}
+\frac{R\xi^{(n-1)}_j\tilde{v}_x}{(\tilde{v}+\phi^{(n-1)}_j)^2}\\
&g_2^{(n-1)}(\phi_j^{(n-1)},\xi_j^{(n-1)},\phi_{jx}^{(n-1)},\psi_{jx}^{(n-1)},\xi_{jx}^{(n-1)})\\
&=-k\frac{\xi^{(n-1)}_{jx}(\tilde{v}_x+\phi^{(n-1)}_{jx})}{(\tilde{v}+\phi^{(n-1)}_j)^2}
-\tilde{p}\psi_{jx}^{(n-1)}-(\tilde{u}_x+\psi_{jx}^{(n-1)})(\frac{R\xi^{(n-1)}_j-\tilde{p}\phi^{(n-1)}_j}
{\tilde{v}+\phi^{(n-1)}_j})+k(\frac{\tilde{\theta}_x}{\tilde{v}+\phi^{(n-1)}_j})_x.
\end{aligned}
\end{equation}
We now assume that $M_0$ suitably small, if $g_2^{(n-1)}\in C(0,t_0;H^{m-1})$, and $\xi_{0j}\in H^{m},$ then there exists a unique local solution $\xi_{j}^{(n)}$ to (\ref{3.8}) satisfying 

\begin{equation}\label{3.11}
\xi^{(n)}_j\in C(0,T_0;H^m)\cap C^1(0,T_0;H^{m-2})\cap L^2(0,T_0;H^{m+1}).
\end{equation}
Making use of this, if $(\phi^{(n-1)}_j(t,x),\psi^{(n-1)}_j(t,x),\xi^{(n-1)}_j(t,x))\in\mathbb{X}_{\frac{3v_+}{8},\frac{3\theta_-}{8},\tilde{C}M},$
from system (\ref{3.8}), we immediately get that

\begin{equation}\label{3.12}
\begin{aligned}
&\|\xi^{(n)}_j(t)\|_2^2+\|\xi^{(n)}_{jt}(t)\|_1^2+
\int_0^{t_0}\|\xi^{(n)}_{jx}(\tau)\|_2^2+\|\xi^{(n)}_{jt}(\tau)\|_2^2d\tau\\
&\leq Ce^{C(v_+,\theta_-,M_0)t_0}
(\|\xi_{0j}\|_2^2+\|\xi_{jt}(0)\|_1^2+C(v_+,\theta_-,M_0)t_0).
\end{aligned}
\end{equation}
Then a direct computation on $(\ref{3.8})_{1,2}$ with (\ref{3.12}) also tell us

\begin{equation}\label{3.13}
\begin{aligned}
&\|(\phi^{(n)}_j,\psi^{(n)}_j)(t)\|_2^2+\|(\phi^{(n)}_{jt},\psi^{(n)}_{jt})(t)\|_1^2\leq Ce^{C(v_+,\theta_-,M_0)t_0}(\|\phi_{0j},\psi_{0j}\|_2^2
+\|(\phi_{jt},\psi_{jt})(0)\|_1^2\\
&+C(v_+,\theta_-,M_0)t_0
+\int_0^{t_0}\|\xi^{(n-1)}_{jxxx}(\tau)\|^2+\|\xi^{(n-1)}_{jtxx}(\tau)\|^2d\tau).
\end{aligned}
\end{equation}
Combing (\ref{3.12}) and (\ref{3.13}), as long as $t_0$ suitably small, we finally get
\begin{equation}\label{3.14}
\|(\phi^{(n)}_j,\psi^{(n)}_j,\xi^{(n)}_j)(t)\|_2+\|(\phi_{jt}^{(n)},\psi_{jt}^{(n)},\xi_{jt}^{(n)})(t)\|_1\leq \tilde{C}M.
\end{equation}
And $\inf_{\mathbb{R}_+\times[0,T]}(\tilde{v}+\phi_j^{(n)})(t,x)\geq \frac{3v_+}{8},
\inf_{\mathbb{R}_+\times[0,T]}(\tilde{\theta}+\xi_j^{(n)})(t,x)\geq \frac{3\theta_-}{8}$.
That is to say the sequence $(\phi^{(n)}_j,\psi^{(n)}_j,\xi^{(n)}_j)$ is uniformly bounded in the function space $\mathbb{X}_{\frac{3v_+}{8},\frac{3\theta_-}{8},\tilde{C}M}(0,t_0).$ By using the same method in (\cite{Kawa}), we can  finally prove that $(\phi^{(n)}_j,\psi^{(n)}_j,\xi^{(n)}_j)$ has a subsequence $(\phi^{(n')}_j,\psi^{(n')}_j,\xi^{(n')}_j)\rightarrow (\phi_j,\psi_j,\xi_j)\in \mathbb{X}_{\frac{3v_+}{8},\frac{3\theta_-}{8},\tilde{C}M}(0,t_0)$ as $n'\rightarrow+\infty$. Again, we let $j\rightarrow\infty,$ we obtain the desired unique local solution $(\phi,\psi,\xi)(t,x)\in
\mathbb{X}_{\frac{3v_+}{8},\frac{3\theta_-}{8},\tilde{C}M}$ as long as $t_0$ so small. Thus Proposition \ref{p1} has been proved.

Note that we could let $M_0(\leq\min\{\frac{5v_+}{8},\frac{5\theta_-}{8}\})$ so small, then
from the proof of  Proposition \ref{p1} and the Sobolev's inequality, the lower bounds of 
$v, \theta$ satisfy following condition obviously.
\begin{equation}\label{3.15-1}
\inf_{\mathbb{R}_+\times[0,t_0]}(\tilde{v}+\phi(t,x))\geq \frac{3v_+}{8},
\quad\quad\inf_{\mathbb{R}_+\times[0,t_0]}(\tilde{\theta}+\xi(t,x))\geq \frac{3\theta_-}{8}.
\end{equation}

\subsection{Stability of the 3-rarefaction wave}
In this section, a priori estimates will be listed, see Proposition \ref{p2}. As long as this Proposition be
proved, then Proposition \ref{p1}, Proposition \ref{p2} and Lemma \ref{L1} imply Theorem \ref{t1}
immediately.

\begin{Proposition}\label{p2} (A priori estimates) For $0<\eps\leq 1$ listed in (\ref{2.6}), if $\eps\leq\eps_1$ for small positive constant $\eps_1$,  there exist  positive constants  $\va_0(\tilde{C}\va_0\leq M_0),$ and $C_0$ which depend on only the initial data such that the following statements hold. If $(\phi,\psi,\xi)(t,x)\in
\mathbb{X}_{\frac{3v_+}{8},\frac{3\theta_-}{8},\va_0}(0,T)$ is a solution to (\ref{3.2}) for some $T>0$, then it holds that

\begin{equation}\label{3.15}
\begin{aligned}
&\|(\phi,\psi,\xi)(t)\|_2^2+\|(\phi_t,\psi_t,\xi_t)(t)\|_1^2+\int_0^t \|(\phi_x,\psi_x)(\tau)\|_1^2+\|\xi_x(\tau)\|_2^2+\|\xi_t(\tau)\|_2^2\,d\tau\\
 &\leq  C_0(\|(\phi_0,\psi_0,\xi_0)\|_2^2+\|(\phi_t,\psi_t,\xi_t)(0)\|_1^2
 +\eps^{\frac{1}{4}}).
\end{aligned}
\end{equation}
\end{Proposition}

Once Proposition \ref{p2} is proved, choosing $\eps,\|(\phi_0,\psi_0,\xi_0)\|_2^2+\|(\phi_t,\psi_t,\xi_t)(0)\|_1^2$ suitably small such that

\begin{equation}\label{3.16}
\begin{aligned}
&\|(\phi_0,\psi_0,\xi_0)\|_2+\|(\phi_t,\psi_t,\xi_t)(0)\|_1\leq\frac{\va_0}{\tilde{C}},\\
&C_0(\|(\phi_0,\psi_0,\xi_0)\|_2^2+\|(\phi_t,\psi_t,\xi_t)(0)\|_1^2+\eps^{\frac{1}{4}})\leq(\frac{\va_0}{\tilde{C}})^2.
\end{aligned}
\end{equation}

Then combing Proposition \ref{p1} and  Proposition \ref{p2}, we can construct the global solution 
$(\phi,\psi,\xi)(t,x)\in\mathbb{X}_{\frac{3v_+}{8},\frac{3\theta_-}{8},\frac{\va_0}{\tilde{C}}}(0,+\infty).$ Moreover, it satisfy the estimate (\ref{3.15}) for all $t>0.$ Then from (\ref{3.15}) and the system (\ref{3.2}), we can get

\begin{equation}\label{3.17}
\begin{aligned}
\int_0^{+\infty}\|(\phi_x,\psi_x,\xi_x)(\tau)\|^2d\tau
+\int_0^{+\infty}\frac{d}{dt}\|(\phi_x,\psi_x,\xi_x)(\tau)\|^2<+\infty,
\end{aligned}
\end{equation}
which  together with the Sobolev's inequality leads to the asymptotic behavior of the solution
\begin{equation}\label{3.18}
\lim_{t\rightarrow+\infty}\sup_{x\in\mathbb{R}_+}|(\phi,\psi,\xi)(t,x)|= 0.
\end{equation}
Combing Lemma \ref{L1}, we finally get Theorem \ref{t1}. Now we mainly proof Proposition \ref{p2}.\\
Proof of Proposition \ref{p2}:\\
First, we show the basic estimates.

\begin{Lemma} Under the same assumptions listed in Proposition \ref{p2}, if $\eps, \va_0$ are suitably small, it holds that 
\begin{equation}\label{3.19}
\begin{aligned}
&\|(\phi,\psi,\xi)(t)\|^2+\IT\|\xi_x(\tau)\|^2
\leq C(v_+,\theta_-)(\|(\phi_0,\psi_0,\zeta_0)\|^2+\eps^{\frac{1}{4}}[\IT\|\phi_x(\tau)\|^2d\tau+1]),
\end{aligned}
\end{equation}
where $C(v_+,\theta_-)$ is a positive constant which depends on $v_+,\theta_-$,  and it is always true without illustration in below.
\end{Lemma}
Proof. Define the energy form
\begin{equation}\label{3.20}
E=R\hat{\theta}\Phi(\frac{v}{\tilde{v}})+\frac{\psi^2}{2}+\frac{R}{\gamma-1}\tilde{\theta}\Phi(\frac{\theta}{\tilde{\theta}}),
\end{equation}
where $\Phi(s)=s-1-\ln s,$ obviously, there exists a positive constant $C(s)$ such that 
\begin{equation}\label{3.21}
C(s)^{-1}s^2\leq\Phi(s)\leq C(s) s^2.
\end{equation}
Then after computation, we get

\begin{equation}\label{3.22}
\begin{aligned}
&E_t+((p-\tilde{p})\psi)_x=k(\frac{\xi\xi_x}{v\theta})_x-k\frac{\tilde{\theta}}{v\theta^2}\xi_{x}^2+H
\end{aligned}
\end{equation}
where
\begin{equation}\label{3.23}
\begin{aligned}
&H=\{-\tilde{p}\tilde{u}_x(\frac{\phi^2}{v\tilde{v}}+\frac{(p-\tilde{p})\xi}{\tilde{p}\theta}-\frac{\xi^2}{\tilde{\theta}\theta}+(\ga-1)\Phi(\frac{v}{\tilde{v}})+\Phi(\frac{\theta}{\tilde{\theta}}))\} \\
&+\{-k\frac{\tilde{\theta}_x\phi_x\xi}{v^2\theta}+k\frac{\tilde{\theta}_x\xi\xi_x}{v\theta^2}+k\frac{\tilde{\theta}_{xx}\xi}{v\theta}
-k\frac{\tilde{\theta}_x\tilde{v}_x\xi}{v^2\theta}\}\\
&=H_1+H_2
\end{aligned}
\end{equation}

The form of $p(v,s)$ in (\ref{2.2-1}) tell us

\begin{equation}\label{3.24}
\begin{aligned}
&0\leq p(v,s)-p(\tilde{v},\tilde{s})-p_v(\tilde{v},\tilde{s})\phi-p_s(\tilde{v},\tilde{s})(s-\tilde{s})\\
&=\tilde{p}(\frac{\phi^2}{v\tilde{v}}+\frac{(p-\tilde{p})\xi}{\tilde{p}\theta}-\frac{\xi^2}{\tilde{\theta}\theta}+(\ga-1)\Phi(\frac{v}{\tilde{v}})
+\Phi(\frac{\theta}{\tilde{\theta}})),
\end{aligned}
\end{equation}
since $\tilde{u}_x\geq0$, 

\begin{equation}\label{3.25}
\int_{0}^{t}\int_{\mathbb{R}_+}H_1dxd\tau\leq 0.
\end{equation}

Besides that,

\begin{equation}\label{3.26}
\begin{aligned}
&\int_{0}^{t}\int_{\mathbb{R}_+}H_2dxd\tau \\
&\leq c(v_+,\theta_-)\eps^{-\frac{1}{4}}\int_{0}^{t}\|\tilde{\theta}_x\|_{L^{\infty}}^2\|\xi\|^2d\tau
+c(v_+,\theta_-)\int_{0}^{t}\|\xi\|_{L^{\infty}}(\|\tilde{u}_x\|^2+\|\tilde{\theta}_{xx}\|_{L^1})d\tau\\
&+\frac{1}{8}k\int_{0}^{t}\int_{\mathbb{R}_+}\frac{\tilde{\theta}}{v\theta^2}\xi_{x}^2dxd\tau+\eps^{\frac{1}{4}}\int_{0}^{t}\|\phi_x(\tau)\|^2d\tau\\
&\leq c(v_+,\theta_-)(\eps^{\frac{1}{4}}\int_0^t(1+\tau)^{-\frac{3}{2}}\|\xi\|^2d\tau+
\int_{0}^{t}\|\xi\|^{\frac{2}{3}}(\|\tilde{u}_x\|^{\frac{8}{3}}+\|\tilde{\theta}_{xx}\|_{L^1}^{\frac{4}{3}})d\tau\\
&+\frac{1}{4}k\int_{0}^{t}\int_{\mathbb{R}_+}\frac{\tilde{\theta}}{v\theta^2}\xi_{x}^2dxd\tau+\eps^{\frac{1}{4}}\int_{0}^{t}\|\phi_x(\tau)\|^2d\tau,\\
\end{aligned}
\end{equation}
where
\begin{equation}\label{3.27}
\begin{aligned}
&\int_{0}^{t}\|\xi\|^{\frac{2}{3}}(\|\tilde{u}_x\|^{\frac{8}{3}}+\|\tilde{\theta}_{xx}\|_{L^1}^{\frac{4}{3}})d\tau
\leq C\eps^{\frac{1}{4}}\int_{0}^{t}(1+\tau)^{(\frac{13}{12})(-1+\frac{1}{q})}(1+\|\xi\|^2)d\tau.\\
\end{aligned}
\end{equation}
Here we use the fact that $|\tilde{\theta}_x|\leq C(\ga-1)|\tilde{u}_x|,|\tilde{v}_x|\leq C|\tilde{u}_x|.$
Combing the results (\ref{3.25})-(\ref{3.27}), choosing $q$ suitably large such that $\frac{13}{12}(-1+\frac{1}{q})>1,$ then
 integrating (\ref{3.22}) over $[0,t]\times\mathbb{R}_+,$ we immediately get that
\begin{equation}\label{3.28}
\begin{aligned}
&\|(\phi,\psi,\xi)(t)\|^2+\int_{0}^{t}\|\xi_x(\tau)\|^2\\
&\leq C(v_+,\theta_-)(\|(\phi_0,\psi_0,\xi_0)\|^2+\eps^{\frac{1}{4}}\int_{0}^{t}\|\phi_x(\tau)\|^2d\tau+\eps^{\frac{1}{4}}\|\xi(t)\|^2+\eps^{\frac{1}{4}}).
\end{aligned}
\end{equation} 
Choosing $\eps^{\frac{1}{4}}$ suiably small, we finally get 
\begin{equation}\label{3.29}
\|(\phi,\psi,\xi)(t)\|^2+\int_{0}^{t}\|\xi_x(\tau)\|^2
\leq C(v_+,\theta_-)(\|(\phi_0,\psi_0,\xi_0)\|^2+\eps^{\frac{1}{4}}\int_{0}^{t}\|\phi_x(\tau)\|^2d\tau+\eps^{\frac{1}{4}}).
\end{equation}
That is (\ref{3.19}) in our Lemma 3.1.

Secondly, we show the one derivative estimates in $L^2.$
\begin{Lemma}Under the same assumptions listed in Proposition \ref{p2}, if $\eps, \va_0$ are suitably small, it holds that
\begin{equation}\label{3.30}
\begin{aligned}
&\|(\phi_x,\psi_x,\xi_x)(t)\|^2+\int_0^t\|\xi_{xx}(\tau)\|^2d\tau
\leq C(v_+,\theta_-)(\|(\phi_0,\psi_0,\xi_0)\|_1^2\\
&+(\va_0+\eps^{\frac{1}{4}}+\eta_1)\int_{0}^{t}\|(\phi_x,\psi_x)(\tau)\|^2d\tau
+(\va_0+\eps^{\frac{1}{4}})\int_0^t\|(\phi_{xx},\psi_{xx})(\tau)\|^2d\tau+\eps^{\frac{1}{4}}). \\
\end{aligned}
\end{equation}
where $\eta_1$ is some small positive constant which will be determined later.
\end{Lemma}
Proof.  Let $(\ref{3.2})_{1x}\times\frac{\tilde{p}}{v}\phi_x, (\ref{3.2})_2\times-\psi_{xx}, (\ref{3.2})_3\times-\frac{\xi_{xx}}{\theta}$, then 

\begin{equation}
\begin{aligned}\label{3.31}
&(\frac{1}{2}\frac{\tilde{p}}{v}\phi_x^2+\frac{1}{2}\psi_x^2+\frac{1}{2}\frac{R}{\gamma-1}
\frac{\xi_x^2}{\theta})_t
-(\psi_t\psi_x+\frac{R}{\gamma-1}\xi_t\frac{\xi_x}{\theta}+\frac{R\xi_x\psi_x}{v})_x\\
&-\frac{1}{2}(\frac{\tilde{p}}{v})_t\phi_x^2
+\frac{R}{\gamma-1}(\frac{1}{\theta})_x\xi_t\xi_x
-\frac{1}{2}\frac{R}{\gamma-1}(\frac{1}{\theta})_t\xi_x^2-\tilde{u}_x(p-\tilde{p})\frac{\xi_{xx}}{\theta}\\
&=-k(\frac{\xi_x}{v})_x \frac{\xi_{xx}}{\theta}-k(\frac{\tilde{\theta}_x}{v})_x\frac{\xi_{xx}}{\theta}+F_1,
\end{aligned}
\end{equation}

where $F_1$ is following:
\begin{equation}\label{3.32}
F_1=-(\frac{R}{v})_x\psi_x\xi_x+(\frac{R}{v})_x\xi\psi_{xx}-(\frac{\tilde{p}}{v})_x\phi\psi_{xx}.
\end{equation}
Integrating (\ref{3.31}) over $[0,t]\times\mathbb{R}_+,$ and by the Sobelov's inequality, 

\begin{equation}\label{3.33}
\begin{aligned}
&\|(\phi_x,\psi_x,\xi_x)(t)\|^2+\int_{0}^{t}\|\xi_{xx}(\tau)\|^2d\tau\\
&\leq C(v_+,\theta_-)(\|(\phi_0,\psi_0,\xi_0)\|_1^2+\int_{0}^{t}|\xi_x\psi_x|(\tau,0)d\tau+
(\va_0+\eps^{\frac{1}{4}})\int_{0}^{t}\|(\phi_x,\psi_x,\xi_x)(\tau)\|_1^2\\
&+\eps^{-\frac{1}{4}}(\int_{0}^{t}\int_{\mathbb{R}_+}|\tilde{u}_x|^2(\phi^2+\xi^2)dxd\tau
+\int_{0}^{t}\|\tilde{\theta}_{xx},\tilde{\theta}_x\tilde{v}_x\|^2d\tau))
\end{aligned}
\end{equation}
Similar as (\ref{3.26}), we immediately get that

\begin{equation}\label{3.34}
\begin{aligned}
&\eps^{-\frac{1}{4}}(\int_{0}^{t}\int_{\mathbb{R}_+}|\tilde{u}_x|^2(\phi^2+\xi^2)dxd\tau+\int_{0}^{t}\|\tilde{\theta}_{xx},\tilde{\theta}_x\tilde{v}_x\|^2d\tau)\\
&\leq C\eps^{\frac{1}{4}}(\int_{0}^{t}(1+\tau)^{-\frac{3}{2}}\|\phi,\xi\|^2d\tau
+\int_{0}^{t}(1+\tau)^{-\frac{5}{2}}+(1+\tau)^{(\frac{5}{3})(-1+\frac{1}{q})}d\tau).
\end{aligned}
\end{equation}

For the boundary term, we deal with it in following way.  On the one hand,

\begin{equation}\label{3.34-1}
\begin{aligned}
&\int_0^t|\xi_x\psi_x|(\tau,0)d\tau \\
&\leq C\int_0^t\|\xi_x(\tau)\|^2d\tau
+\frac{1}{4}\int_{0}^{t}\|\xi_{xx}(\tau)\|^2d\tau
+\eta_1\int_{0}^{t}|\psi_x|^2(\tau,0)d\tau,
\end{aligned}
\end{equation}
and  $\eta_1$ is some small constant which will be determined later, on the other hand,  $(\ref{3.2})_{1x}\times\psi_{x}+(\ref{3.2})_{2x}\times\phi_x$ tell us
\begin{equation}\label{3.34-2}
(\phi_x\psi_x)_t-(\frac{1}{2}\psi_x^2+\frac{\tilde{p}}{2v}\phi_x^2)_x=F_2,
\end{equation}
where $F_2$ is 
\begin{equation}\label{3.34-3}
F_2=2\frac{R\xi_xv_x}{v^2}\phi_x+\frac{3}{2}(\frac{\tilde{p}}{v})_x\phi_x^2+(\frac{\tilde{p}}{v})_{xx}\phi\phi_x-\frac{R\xi_{xx}}{v}\phi_x
-(\frac{R}{v})_{xx}\xi\phi_x.
\end{equation}
Integrating (\ref{3.34-2}) over $[0,t]\times\mathbb{R}_+$, then

\begin{equation}\label{3.34-4}
\begin{aligned}
&\int_{0}^{t}\phi_x^2(\tau,0)+\psi_x^2(\tau,0)d\tau\\
&\leq  C(v_+,\theta_-)(\|(\phi_0,\psi_0)\|_1^2+\|(\phi_x,\psi_x)(t)\|^2+(\va_0
+\eps^{\frac{1}{4}})\int_{0}^{t}\|(\xi_x,\phi_{xx})(\tau)\|^2\\
&+\int_0^t\|(\phi_x,\xi_{xx})(\tau)\|^2d\tau)+\eps^{-\frac{1}{4}}\int_{0}^{t}
(\|\tilde{u}_{xx}\|_{L^\infty}^2+\|\tilde{u}_x\|^4_{L^\infty})\|\phi,\xi\|^2d\tau
\end{aligned}
\end{equation}
Simiar as before, we have 
\begin{equation}\label{3.34-5}
\begin{aligned}
&\eps^{-\frac{1}{4}}\int_{0}^{t}
(\|\tilde{u}_{xx}\|_{L^\infty}^2+\|\tilde{u}_x\|^4_{L^\infty})\|\phi,\xi\|^2d\tau
\leq c\eps^{\frac{1}{4}}\int_{0}^{t}(1+\tau)^{-\frac{3}{2}}\|\phi,\xi\|^2d\tau
\end{aligned}
\end{equation}

Substituting (\ref{3.19}), (\ref{3.34})-(\ref{3.34-1}) and (\ref{3.34-4})-(\ref{3.34-5}) into (\ref{3.33}),  choosing $\eps, \va_0$ suitably small,  $q$ suitably large, we have

\begin{equation}\label{3.35}
\begin{aligned}
&\|(\phi_x,\psi_x,\xi_x)(t)\|^2+\int_0^t\|\xi_{xx}(\tau)\|^2d\tau\\
&\leq C(v_+,\theta_-)(\|(\phi_0,\psi_0,\xi_0)\|_1^2+\eta_1\|(\phi_x,\psi_x)(t)\|^2
+(\va_0+\eps^{\frac{1}{4}})\int_{0}^{t}\|(\phi_x,\psi_x)(\tau)\|_1^2d\tau+\eps^{\frac{1}{4}}\\
&+\eta_1\int_0^t\|(\phi_x,\xi_{xx})(\tau)\|^2d\tau).
\end{aligned}
\end{equation}
Again let $\eta_1$ suitably small, we get
\begin{equation}\label{3.36}
\begin{aligned}
&\|(\phi_x,\psi_x,\xi_x)(t)\|^2+\int_0^t\|\xi_{xx}(\tau)\|^2d\tau
\leq C(v_+,\theta_-)(\|(\phi_0,\psi_0,\xi_0)\|_1^2\\
&+(\va_0+\eps^{\frac{1}{4}}+\eta_1)\int_{0}^{t}\|(\phi_x,\psi_x)(\tau)\|^2d\tau
+(\va_0+\eps^{\frac{1}{4}})\int_0^t\|(\phi_{xx},\psi_{xx})(\tau)\|^2d\tau+\eps^{\frac{1}{4}}). \\
\end{aligned}
\end{equation}
This is our estimate $(\ref{3.30}).$


\begin{Lemma}Under the same assumptions listed in Lemma 3.2, if $\eps, \va_0,\eta_1$ are suitably small, it holds that
	
\begin{equation}\label{3.37}
\begin{aligned}
&\int_0^t\|(\phi_x,\psi_x)(\tau)\|^2d\tau \\
&\leq C(v_+,\theta_-)(\|(\phi_0,\psi_0,\xi_0)\|_1^2+(\va_0+\eps^{\frac{1}{4}})\int_{0}^{t}\|(\phi_{xx},\psi_{xx})(\tau)\|^2d\tau+\eps^{\frac{1}{4}}),
\end{aligned}
\end{equation}
\end{Lemma}
Proof.  $(\ref{3.2})_2\times-\frac{1}{2}\phi_x, (\ref{3.2})_3\times\frac{\psi_x}{p}$ imply that 

\begin{equation}\label{3.38}
\begin{aligned}
&\frac{\tilde{p}\phi_x^2}{2v}+\frac{1}{2}[(\frac{\tilde{p}}{v})_x\phi\phi_x-(\frac{R\xi}{v})_x\phi_x]=\frac{1}{2}[(\psi\phi_x)_t-(\psi\phi_t)_x+\psi_x^2] \\
&(\frac{R}{\ga-1}\frac{\psi_x}{p}\xi)_t-\frac{R}{\ga-1}(\frac{1}{p})_t\xi\psi_x
-(\frac{R}{\ga-1}\frac{\psi_t\xi}{p})_x+\frac{R}{\ga-1}(\frac{\xi}{p})_x\psi_t+\psi_x^2\\
&+\frac{\tilde{u}_x(p-\tilde{p})}{p}\psi_x=k(\frac{\theta_x}{v})_x\frac{\psi_x}{p}
\end{aligned}
\end{equation}

Integrating $(\ref{3.38})_1+(\ref{3.38})_2$ over $[0,t]\times\mathbb{R}_+,$ and by Sobelov's
inequality and previous results  $(\ref{3.19}), (\ref{3.30})$, 

\begin{equation}\label{3.39}
\begin{aligned}
&\int_0^t\|(\phi_x,\psi_x)(\tau)\|^2d\tau 
\leq C(v_+,\theta_-)(\|(\phi_0,\psi_0,\xi_0)\|_1^2\\
&+(\va_0+\eps^{\frac{1}{4}}+\eta_1)\int_{0}^{t}\|(\phi_x,\psi_x)(\tau)\|^2d\tau+(\va_0+\eps^{\frac{1}{4}})\int_0^t\|(\phi_{xx},
\psi_{xx})(\tau)\|^2+\eps^{\frac{1}{4}}).
\end{aligned}
\end{equation}
 Again choosing $\va_0,\eps,\eta_1$ suitably small, we finally get 

\begin{equation}\label{3.41}
\begin{aligned}
&\int_0^t\|(\phi_x,\psi_x)(\tau)\|^2d\tau \\
&\leq C(v_+,\theta_-)(\|(\phi_0,\psi_0,\xi_0)\|_1^2+(\va_0+\eps^{\frac{1}{4}})\int_{0}^{t}\|(\phi_{xx},\psi_{xx})(\tau)\|^2d\tau+\eps^{\frac{1}{4}}),
\end{aligned}
\end{equation}
which is our estimate (\ref{3.37}).
Combining Lemma 3.1-Lemma 3.3, our estimates are listed as follows

\begin{equation}\label{3.42}
\begin{aligned}
&\|(\phi,\psi,\xi)(t)\|_1^2+\int_{0}^{t}\|(\phi_x,\psi_x)(\tau)\|^2+\|\xi_x(\tau)\|_1^2d\tau \\
&\leq C(v_+,\theta_-)(\|(\phi_0,\psi_0,\xi_0)\|_1^2+(\va_0
+\eps^{\frac{1}{4}})\int_{0}^{t}\|(\phi_{xx},\psi_{xx})(\tau)\|^2d\tau+\eps^{\frac{1}{4}} )
\end{aligned}
\end{equation}

By our system (\ref{3.2}),  from (\ref{3.42}) and previous computations, we can deduce that

\begin{equation}\label{3.43}
\begin{aligned}
&\int_0^t\|(\phi_t,\psi_t,\xi_t)(\tau)\|^2d\tau\\
&\leq C(v_+,\theta_-)(\|(\phi_0,\psi_0,\xi_0)\|_1^2+(\va_0
+\eps^{\frac{1}{4}})\int_{0}^{t}\|(\phi_{xx},\psi_{xx})(\tau)\|^2d\tau+\eps^{\frac{1}{4}} )
\end{aligned}
\end{equation}

\begin{Lemma}Under the same assumptions listed in Proposition 3.2, if $\eps, \va_0$ are suitably small, it holds that
	
\begin{equation}\label{3.44}
\begin{aligned}
&\|(\phi_t,\psi_t,\xi_t)(t)\|^2+\int_{0}^{t}\|\xi_{tx}(\tau)\|^2d\tau\\
&\leq C(v_+,\theta_-)(\|(\phi_0,\psi_0,\xi_0)\|_1^2+\|(\phi_t,\psi_t,\xi_t)(0)\|^2+(\va_0+\eps^{\frac{1}{4}})\int_{0}^{t}
\|(\phi_{xx},\psi_{xx})(\tau)\|^2d\tau+\eps^{\frac{1}{4}})
\end{aligned}
\end{equation}
\end{Lemma}
Proof.  $(\ref{3.2})_{1t}\times\frac{\tilde{p}}{v}\phi_t,(\ref{3.2})_{2t}\times\psi_t,
(\ref{3.2})_{3t}\times\frac{\xi_t}{\theta}$

\begin{equation}\label{3.45}
\begin{aligned}
&(\frac{1}{2}\frac{\tilde{p}}{v}\phi_t^2+\frac{1}{2}\psi_t^2+\frac{1}{2}\frac{R}{\gamma-1}\frac{\xi_t^2}{\theta})_t+(\frac{R\xi_t-\tilde{p}\phi_t}{v}\psi_t)_x
+(\tilde{u}_x(p-\tilde{p}))_t\frac{\xi_t}{\theta}\\
&=k((\frac{\theta_x}{v})_t\frac{\xi_t}{\theta})_x-k(\frac{\theta_x}{v})_t\frac{\xi_{tx}}{\theta}+F_3,
\end{aligned}
\end{equation}
where $F_3$ is 

\begin{equation}\label{3.46}
\begin{aligned}
&F_3=\frac{1}{2}(\frac{\tilde{p}}{v})_t\phi_t^2+\frac{1}{2}\frac{R}{\ga-1}(\frac{1}{\theta})_t\xi_t^2+(\frac{\tilde{p}}{v})_t\phi_x\psi_t-(\frac{R}{v})_t\xi_x\psi_t\\
&+[(\frac{\tilde{p}}{v})_{tx}\phi-(\frac{R}{v})_{tx}\xi]\psi_t
-p_t\psi_x\frac{\xi_t}{\theta}-k(\frac{\theta_x}{v})_t(\frac{1}{\theta})_x\xi_t.
\end{aligned}
\end{equation}
Integrating (\ref{3.45}) over $[0,t]\times\mathbb{R}_+$ and note that $\psi_t(\tau,0)=0, \xi_t(\tau,0)=0,$ then  making use of (\ref{3.42}), (\ref{3.43})  and previous results, when $\va_0, \eps$  are suitably small,  

\begin{equation}\label{3.47}
\begin{aligned}
&\|(\phi_t,\psi_t,\xi_t)(t)\|^2+\int_{0}^{t}\|\xi_{tx}(\tau)\|^2d\tau\\
&\leq C(v_+,\theta_-)(\|(\phi_t,\psi_t,\xi_t)(0)\|^2+\|(\phi_0,\psi_0,\xi_0)\|_1^2
+(\va_0+\eps^{\frac{1}{4}})\int_{0}^{t}\|(\phi_{xx},\psi_{xx})(\tau)\|^2d\tau+\eps^{\frac{1}{4}})
\end{aligned}
\end{equation}
That is $(\ref{3.44}).$

Lastly, we show  the higher order derivative estimates.

\begin{Lemma}Under the same assumptions listed in Proposition 3.2, if $\eps, \va_0$ are suitably small, it holds that
	
\begin{equation}\label{3.48}
\begin{aligned}
&\|(\phi_{xx},\psi_{xx},\xi_{xx})(t)\|^2+\int_{0}^{t}\|\xi_{xxx}(\tau)\|^2d\tau\\
&\leq C(v_+,\theta_-)(\|(\phi_0,\psi_0,\xi_0)\|_2^2+\eps^{\frac{1}{4}}+
(\eta_2+\eta_3+\va_0+\eps^{\frac{1}{4}})\int_{0}^{t}\|(\phi_{xx},\psi_{xx})(\tau)\|^2d\tau\\
&+\eta_3^{-1}\int_{0}^{t}\xi_{tx}^2(\tau,0)d\tau).
\end{aligned}
\end{equation}
where $\eta_2, \eta_3$ are small positive constants which will be determined later.
\end{Lemma}
Proof. $(\ref{3.2})_{1x}\times-\frac{\tilde{p}}{v}\phi_{xxx}, (\ref{3.2})_{2xx}\times\psi_{xx},
(\ref{3.2})_{3x}\times-\frac{\xi_{xxx}}{\theta}$ tell us 

\begin{equation}\label{3.49}
\begin{aligned}
&\frac{1}{2}(\frac{\tilde{p}}{v}\phi_{xx}^2+\psi_{xx}^2+\frac{R}{\ga-1}\frac{1}{\theta}\xi_{xx}^2)_t-(\frac{\tilde{p}}{v}\phi_{tx}\phi_{xx}+\frac{R}{\ga-1}\xi_{tx}\frac{\xi_{xx}}{\theta})_x
-(\tilde{u}_x(p-\tilde{p}))_x\frac{\xi_{xxx}}{\theta}\\
&+k\frac{\xi_{xxx}^2}{v\theta}=F_4-2k(\frac{1}{v})_{x}\xi_{xx}\frac{\xi_{xxx}}{\theta}-k(\frac{1}{v})_{xx}\xi_x\frac{\xi_{xxx}}{\theta}-k(\frac{\tilde{\theta}_x}{v})_{xx}\frac{\xi_{xxx}}{\theta},
\end{aligned}
\end{equation}
where $F_4$ is 
\begin{equation}\label{3.50}
\begin{aligned}
&F_4=\frac{1}{2}(\frac{\tilde{p}}{v})_t\phi_{xx}^2+\frac{1}{2}(\frac{R}{\ga-1}\frac{1}{\theta})_t\xi_{xx}^2
-(\frac{\tilde{p}}{v})_x\phi_{tx}\phi_{xx}-\frac{R}{\ga-1}(\frac{1}{\theta})_x\xi_{tx}\xi_{xx}+p_x\psi_x\frac{\xi_{xxx}}{\theta}\\
&+3((\frac{\tilde{p}}{v})_{xx}\phi_x+(\frac{\tilde{p}}{v})_x\phi_{xx})\psi_{xx}+(\frac{\tilde{p}}{v})_{xxx}\phi\psi_{xx}-3((\frac{R}{v})_{xx}\xi_x+(\frac{R}{v})_x\xi_{xx})\psi_{xx}-(\frac{R}{v})_{xxx}\xi\psi_{xx}.
\end{aligned}
\end{equation}
Integrating (\ref{3.49}) over $[0,t]\times\mathbb{R}_+,$ by $(\ref{3.42})$ and $(\ref{3.44})$

\begin{equation}\label{3.51}
\begin{aligned}
&\|(\phi_{xx},\psi_{xx},\xi_{xx})(t)\|^2+\int_{0}^{t}\|\xi_{xxx}(\tau)\|^2d\tau\\
&\leq C(v_+,\theta_-)(\|(\phi_0,\psi_0,\xi_0)\|_2^2+|\int_{0}^{t}H(\tau,0)d\tau|\\
&+(\va_0+\eps^{\frac{1}{4}})\int_{0}^{t}\|(\phi_{xx},\psi_{xx})(\tau)\|^2+\|\xi_{xxx}(\tau)\|^2d\tau+\eps^{\frac{1}{4}}\\
&+\int_{0}^{t}\int_{\mathbb{R}_+}((\frac{\tilde{p}}{v})_{xxx}\phi-(\frac{R}{v})_{xxx}\xi)\psi_{xx}dxd\tau) .
\end{aligned}
\end{equation}
where 
\begin{equation*}
H(\tau,0)=\frac{R}{\ga-1}\xi_{tx}\xi_{xx}(\tau,0)+\frac{\tilde{p}}{v}\phi_{tx}\phi_{xx}(\tau,0)=H_1(\tau,0)+H_2(\tau,0)
\end{equation*}

On the one hand, using the relationship $\phi_{txx}=\psi_{xxx}$ and previous results,  we can get

\begin{equation}\label{3.52}
\begin{aligned}
&\int_{0}^{t}\int_{\mathbb{R}_+}((\frac{\tilde{p}}{v})_{xxx}\phi-(\frac{R}{v})_{xxx}\xi)\psi_{xx}dxd\tau\\ 
&= \int_{0}^{t}\int_{\mathbb{R}_+}[(\frac{\tilde{p}_{xxx}}{v}+3(\frac{1}{v})_x\tilde{p}_{xx}+3(\frac{1}{v})_{xx}\tilde{p}_x)\phi
+(\tilde{p}\phi-R\xi)(\frac{1}{v})_{xxx}]\psi_{xx}\\
&=\int_{0}^{t}\int_{\mathbb{R}_+}[(\frac{\tilde{p}_{xxx}}{v}+3(\frac{1}{v})_x\tilde{p}_{xx}+3(\frac{1}{v})_{xx}\tilde{p}_x)\phi
+(\tilde{p}\phi-R\xi)(-\frac{\tilde{v}_{xxx}}{v^2}+6\frac{v_xv_{xx}}{v^3}-6\frac{v_x^3}{v^4})]\psi_{xx}\\
&-\int_{0}^{t}\int_{\mathbb{R}_+}(\tilde{p}\phi-R\xi)(\frac{\phi_{xxx}}{v^2})\psi_{xx}\\
&\leq C(v_+,\theta_-)((\va_0+\eps^{\frac{1}{4}})\int_{0}^{t}\|(\phi_x,\psi_x)(\tau)
\|_1^2+\eps^{-\frac{1}{4}}\int_{0}^{t}\|\tilde{u}_x\|_{L^\infty}^2\|\phi,\xi\|^2d\tau)\\
&-\int_{0}^{t}\int_{\mathbb{R}_+}((\tilde{p}\phi-R\xi)\frac{\phi_{xx}}{v^2}\psi_{xx})_x
+\int_{0}^{t}\int_{\mathbb{R}_+}(\frac{\tilde{p}\phi-R\xi}{v^2})_x\phi_{xx}\psi_{xx}dxd\tau\\
&+\int_{0}^{t}\int_{\mathbb{R}_+}((\frac{\tilde{p}\phi-R\xi}{2v^2})\phi_{xx}^2)_t-(\frac{\tilde{p}\phi-R\xi}{2v^2})_t\phi_{xx}^2dxd\tau\\
&\leq C(v_+,\theta_-)(\|\phi_0,\psi_0,\xi_0\|_2^2+\va_0\|\phi_{xx}\|^2+\eps^{\frac{1}{4}}
+(\va_0+\eps^{\frac{1}{4}})\int_{0}^{t}\|(\phi_x,\psi_x)(\tau)\|_1^2+\|\xi_{xxx}(\tau)\|^2d\tau\\
&+\va_0\int_{0}^{t}|\phi_{xx}|^2(\tau,0)d\tau),
\end{aligned}
\end{equation}
on the other hand,  we turn to deal with the higher order boundary terms $H(\tau,0)$ in (\ref{3.51}). By the conductivity of the equation $(\ref{3.2})_3$, we know that
\begin{equation}\label{3.53}
\begin{aligned}
\xi_{xx}(\tau,0)=(\frac{1}{k}(R\theta_-\psi_x-\tilde{u}_x\tilde{p}\phi)+\frac{\xi_xv_x}{v}-\tilde{\theta}_{xx}+\frac{\tilde{\theta}_xv_x}{v})(\tau,0).
\end{aligned}
\end{equation}
Hence $H_1(\tau,0)$ can be estimated in following way

\begin{equation}\label{3.54}
\begin{aligned}
&|\int_0^t\frac{R}{\ga-1}\xi_{tx}\xi_{xx}(\tau,0)d\tau| \\
&=|\int_0^t\frac{R}{(\ga-1)}[(\frac{1}{k}(R\theta_-\psi_x-\tilde{u}_x\tilde{p}\phi)+\frac{\xi_xv_x}{v}-\tilde{\theta}_{xx}+\frac{\tilde{\theta}_xv_x}{v})(\tau,0)
\xi_x(\tau,0)]_td\tau\\
&-\int_0^t\frac{R}{\ga-1}(\frac{1}{k}(R\theta_-\psi_x-\tilde{u}_x\tilde{p}\phi)+\frac{\xi_xv_x}{v}-\tilde{\theta}_{xx}+\frac{\tilde{\theta}_xv_x}{v})_t(\tau,0)\xi_x(\tau,0)d\tau|\\
&\leq C(v_+,\theta_-) (\|\phi_0,\psi_0,\xi_0\|_2^2+\|\phi,\psi,\xi\|_1^2+(\eta_2+\va_0+\eps^{\frac{1}{4}})\|\phi_{xx},\psi_{xx},\xi_{xx}\|^2+\eps^{\frac{1}{4}}\\
&+(\va_0+\eps^{\frac{1}{4}})\int_{0}^{t}\|\phi_{x}(\tau),\psi_x(\tau)\|_1^2d\tau+\frac{1}{4}\int_{0}^{t}\|\xi_{xxx}\|^2d\tau
+\int_0^t\|\xi_{x}(\tau)\|_1^2d\tau\\
&+(\eta_2+\va_0+\eps^{\frac{1}{4}})\int_0^t(\phi_{xx}^2+\psi_{xx}^2)(\tau,0)d\tau).
\end{aligned}
\end{equation}
Here $\eta_2$ is a small positive constant which will be determined below.
For the other higher order boundary term  $H_2(\tau,0)$ in $(\ref{3.51})$  which is one of the main difficulties in our estimates, we deal with it in following procedure.  By $(\ref{3.2})_2$, 

\begin{equation}\label{3.55}
\begin{aligned}
\psi_{tt}(\tau,0)+(\frac{R\xi-\tilde{p}\phi}{v})_{tx}(\tau,0)=0,
\end{aligned}
\end{equation}
and $\psi_{tt}(\tau,0)=0$ immediately gives us that

\begin{equation}\label{3.56}
\begin{aligned}
\frac{\tilde{p}}{v}\phi_{tx}(\tau,0)=\frac{R}{v}\xi_{tx}(\tau,0)-((\frac{\tilde{p}}{v})_{tx}\phi+(\frac{\tilde{p}}{v})_t\phi_x+(\frac{\tilde{p}}{v})_x\phi_t
-(\frac{R}{v})_{t}\xi_x)(\tau,0).
\end{aligned}
\end{equation}
Using this, and by Young-inequality, it holds that

\begin{equation}\label{3.57}
\begin{aligned}
&|\int_0^t\frac{\tilde{p}}{v}\phi_{tx}\phi_{xx}(\tau,0)d\tau|\\
&\leq C(v_+,\theta_-)(\|\phi_0,\psi_0,\xi_0\|_2^2+\eps^{\frac{1}{4}}+(\va_0+\eps^{\frac{1}{4}})
\int_0^t(\phi_x^2+\psi_x^2+\xi_x^2+\phi_{xx}^2+\psi_{xx}^2)(\tau,0)d\tau\\
&+(\va_0+\eps^{\frac{1}{4}})\int_{0}^{t}\|(\phi_{xx},\psi_{xx})(\tau)\|^2d\tau+\eta_3\int_{0}^{t}\phi_{xx}^2(\tau,0)d\tau
+\eta_3^{-1}\int_0^t\xi_{tx}^2(\tau,0)d\tau,
\end{aligned}
\end{equation}
where $\eta_3$ is some small positive constant to be determined. Substituting $(\ref{3.52}),(\ref{3.54})$, $(\ref{3.57})$  into $(\ref{3.51})$ and using (\ref{3.42}), we have 

\begin{equation}\label{3.58}
\begin{aligned}
&\|(\phi_{xx},\psi_{xx},\xi_{xx})(t)\|^2+\int_{0}^{t}\|\xi_{xxx}(\tau)\|^2d\tau\\
&\leq C(v_+,\theta_-)(\|(\phi_0,\psi_0,\xi_0)\|_2^2+(\eta_2+\va_0+\eps^{\frac{1}{4}})\|(\phi_{xx},\psi_{xx},\xi_{xx})(t)\|^2\\
&+\eps^{\frac{1}{4}}
+(\va_0+\eps^{\frac{1}{4}})\int_{0}^{t}\|(\phi_{xx},\psi_{xx})(\tau)\|^2d\tau\\
&+(\eta_2+\eta_3+\va_0+\eps^{\frac{1}{4}})\int_{0}^{t}(\phi_{xx}^2+\psi_{xx}^2)(\tau,0)d\tau+\eta_3^{-1}\int_{0}^{t}\xi_{tx}^2(\tau,0)d\tau).
\end{aligned}
\end{equation}

Similar as procedures  (\ref{3.34-2})-(\ref{3.34-4}), let $(\ref{3.2})_{1xx}\times\psi_{xx}+(\ref{3.2})_{2xx}\times\phi_{xx}$, we have 
\begin{equation}\label{3.59}
\begin{aligned}
(\phi_{xx}\psi_{xx})_t-\frac{1}{2}(\frac{\tilde{p}}{v}\phi_{xx}^2+\psi_{xx}^2)_x=F_5,
\end{aligned}
\end{equation} 
where $F_5$ is 
\begin{equation}\label{3.60}
\begin{aligned}
&F_5= \frac{5}{2}(\frac{\tilde{p}}{v})_x\phi_{xx}^2+3(\frac{\tilde{p}}{v})_{xx}\phi_x\phi_{xx}-3((\frac{R}{v})_{xx}\xi_x-\frac{R}{v}\xi_{xxx}\phi_{xx}.\\
&+(\frac{R}{v})_x\xi_{xx})\phi_{xx}+\{((\frac{\tilde{p}}{v})_{xxx}\phi-(\frac{R}{v})_{xxx}\xi)\phi_{xx}\}
\end{aligned}
\end{equation}
Integating $(\ref{3.59})$ over $[0,t]\times\mathbb{R}_+,$ we can get

\begin{equation}\label{3.61}
\begin{aligned}
&\int_{0}^{t}\phi_{xx}^2(\tau,0)+\psi_{xx}^2(\tau,0) d\tau
\leq C(v_+,\theta_-)(\|(\phi_0,\psi_0)\|_2^2+\|(\phi_{xx},\psi_{xx})(t)\|^2+\eps^{\frac{1}{4}}\\
&+(\va_0+\eps^{\frac{1}{4}})\int_0^t\|(\psi_{xx},\xi_{xx})(\tau)\|^2d\tau
+\int_{0}^{t}\|(\phi_{xx},\xi_{xxx})(\tau)\|^2d\tau).
\end{aligned}
\end{equation}
Here the estimate of the last term  in (\ref{3.60}) is similar as $(\ref{3.52}).$ Using
(\ref{3.61}), (\ref{3.58}) turns to

\begin{equation}\label{3.62}
\begin{aligned}
&\|(\phi_{xx},\psi_{xx},\xi_{xx})(t)\|^2+\int_{0}^{t}\|\xi_{xxx}(\tau)\|^2d\tau\\
&\leq C(v_+,\theta_-)(\|(\phi_0,\psi_0,\xi_0)\|_2^2+(\eta_2+\eta_3+\va_0+\eps^{\frac{1}{4}})\|(\phi_{xx},\psi_{xx},\xi_{xx})(t)\|^2+\eps^{\frac{1}{4}}\\
&
+(\eta_2+\eta_3+\va_0+\eps^{\frac{1}{4}})\int_{0}^{t}\|(\phi_{xx},\psi_{xx},\xi_{xxx})(\tau)\|^2d\tau+\eta_3^{-1}\int_{0}^{t}\xi_{tx}^2(\tau,0)d\tau).
\end{aligned}
\end{equation}
Choosing $\eta_2,\eta_3,\va_0,\eps$ suitably small such that
\begin{equation}\label{3.63}
C(v_+,\theta_-)(\eta_2+\eta_3+\va_0+\eps^{\frac{1}{4}})\leq \frac{1}{2},
\end{equation}
we get
\begin{equation}\label{3.66}
\begin{aligned}
&\|(\phi_{xx},\psi_{xx},\xi_{xx})(t)\|^2+\int_{0}^{t}\|\xi_{xxx}(\tau)\|^2d\tau\\
&\leq C(v_+,\theta_-)(\|(\phi_0,\psi_0,\xi_0)\|_2^2+\eps^{\frac{1}{4}}+
(\eta_2+\eta_3+\va_0+\eps^{\frac{1}{4}})\int_{0}^{t}\|(\phi_{xx},\psi_{xx})(\tau)\|^2d\tau\\
&+\eta_3^{-1}\int_{0}^{t}\xi_{tx}^2(\tau,0)d\tau).
\end{aligned}
\end{equation}
That is our estimate $(\ref{3.48}).$

\begin{Lemma}Under the same assumptions listed in Proposition 3.2, if $\eps, \va_0$ are suitably small, it holds that
	
\begin{equation}\label{3.67}
\begin{aligned}
&\int_0^t\|(\phi_{xx},\psi_{xx})(\tau)\|^2d\tau\\
&\leq C(v_+,\theta_-)(\|(\phi_0,\psi_0,\xi_0)\|_2^2+\|(\phi_t,\psi_t,\xi_t)(0)\|^2+\eps^{\frac{1}{4}}+\eta_4\int_{0}^{t}\|\xi_{txx}(\tau)\|^2d\tau).
\end{aligned}
\end{equation}
where $\eta_4$ is a small positive constant which will be determined later.
\end{Lemma}

Proof. Similar as Lemma 3.3, $(\ref{3.2})_{2x}\times-\frac{\phi_{xx}}{2},(\ref{3.2})_{3x}\times\frac{\psi_{xx}}{p}$ give us
\begin{equation}\label{3.68}
\begin{aligned}
&\frac{\tilde{p}}{2v}\phi_{xx}^2+(\frac{\tilde{p}}{v})_x\phi_x\phi_{xx}+(\frac{\tilde{p}}{2v})_{xx}\phi\phi_{xx}
=\frac{1}{2}((\psi_x\phi_{xx})_t-(\psi_x\phi_{tx})_x+\psi_{xx}^2+(\frac{R\xi}{v})_{xx}\phi_{xx})\\
&(\frac{R}{\ga-1}\frac{\psi_{xx}}{p}\xi_x)_t-\frac{R}{\ga-1}(\frac{1}{p})_t\xi_x\psi_{xx}
-(\frac{R}{\ga-1}\frac{\psi_{tx}\xi_x}{p})_x+\frac{R}{\ga-1}(\frac{\xi_x}{p})_x\psi_{tx}+\psi_{xx}^2+p_x\psi_x\frac{\psi_{xx}}{p}\\
&+(\tilde{u}_x(p-\tilde{p}))_x\frac{\psi_{xx}}{p}=k(\frac{\theta_x}{v})_{xx}\frac{\psi_{xx}}{p}.
\end{aligned}
\end{equation}
Then integrating $(\ref{3.68})_1+(\ref{3.68})_2$ over $[0,t]\times\mathbb{R}_+,$ we get 

\begin{equation}\label{3.69}
\begin{aligned}
&\int_0^t\|(\phi_{xx},\psi_{xx})(\tau)\|^2d\tau \leq C(v_+,\theta_-)(\|(\phi_0,\psi_0,\xi_0)\|_2^2+\|(\phi_{x},\psi_{x},\xi_{x})(t)\|_1^2
+\int_{0}^{t}\|\xi_x(\tau)\|_2^2d\tau\\
&+\eps^{\frac{1}{4}}
+(\va_0+\eps^{\frac{1}{4}})\int_0^t\|(\phi_{xx},\psi_{xx})(\tau)\|^2d\tau)+|\int_0^t(\psi_x\phi_{tx})(\tau,0)d\tau|+
|\int_0^t(\frac{\psi_{tx}\xi_x}{p})(\tau,0)d\tau|
\end{aligned}
\end{equation}
For the boundary terms in $(\ref{3.69}),$ making use of $(\ref{3.56}), (\ref{3.34-4}), (\ref{3.42})$ and (\ref{3.48}), (\ref{3.61}) we have 
\begin{equation}\label{3.70}
\begin{aligned}
&|\int_0^t(\psi_x\phi_{tx})(\tau,0)d\tau|\leq C(v_+,\theta_-)(\|\phi_0,\psi_0,\xi_0\|_2^2+ \eps^{\frac{1}{4}}+(\va_0+\eps^{\frac{1}{4}})\int_0^t\|(\phi_{xx},\psi_{xx})(\tau)\|^2\\
&+[(\va_0+\eps^{\frac{1}{4}})\eta_3^{-1}+1]\int_0^t\xi_{tx}^2(\tau,0)d\tau)\\
\end{aligned}
\end{equation}
and

\begin{equation}\label{3.71}
\begin{aligned}
&|\int_0^t(\frac{\psi_{tx}\xi_x}{p})(\tau,0)d\tau|= |\int_{0}^{t}(\frac{\psi_x\xi_x}{p})_t(\tau,0)-\int_0^t((\frac{\xi_x}{p})_t\psi_x)(\tau,0)|d\tau\\
&\leq C(v_+,\theta_-)(\|\phi_0,\psi_0,\xi_0\|_2^2+\|\psi_x,\xi_x\|_1^2+\eps^{\frac{1}{4}}+(\va_0+\eps^{\frac{1}{4}})\int_0^t\|(\phi_{xx},\psi_{xx})(\tau)\|^2d\tau\\
&+\int_{0}^{t}\xi_{tx}^2(\tau,0)d\tau).
\end{aligned}
\end{equation}
Substituting $(\ref{3.70})$ and $(\ref{3.71})$ into $(\ref{3.69}),$ and using $(\ref{3.42})$, $(\ref{3.48})$, we have

\begin{equation}\label{3.72}
\begin{aligned}
&\int_0^t\|(\phi_{xx},\psi_{xx})(\tau)\|^2d\tau \leq C(v_+,\theta_-)(\|(\phi_0,\psi_0,\xi_0)\|_2^2
+\eps^{\frac{1}{4}}\\
&+(\eta_2+\eta_3+\va_0+\eps^{\frac{1}{4}})\int_0^t\|(\phi_{xx},\psi_{xx})(\tau)\|^2d\tau+\eta_3^{-1}\int_0^t\xi_{tx}^2(\tau,0)d\tau).
\end{aligned}
\end{equation}
Again choosing $\eta_2,\eta_3,\va_0,\eps$ suitably small such that
\begin{equation}\label{3.73}
C(v_+,\theta_-)(\eta_2+\eta_3+\va_0+\eps^{\frac{1}{4}})\leq \frac{3}{4},\quad 
C(v_+,\theta_-)\eta_3^{-1}\leq \bar{C}(v_+,\theta_-),
\end{equation}
for some positive constant $\bar{C}(v_+,\theta_-)$ which is larger than $C(v_+,\theta_-),$ the condition (\ref{3.73}) we proposed could be satisfied obviously. Then using (\ref{3.44}), by Sobolev's inequality and Young inequality, (\ref{3.72}) becomes

\begin{equation}\label{3.74}
\begin{aligned}
&\int_0^t\|(\phi_{xx},\psi_{xx})(\tau)\|^2d\tau\\ 
&\leq C_1(v_+,\theta_-)(\|\phi_0,\psi_0,\xi_0\|_2^2+\|(\phi_t,\psi_t,\xi_t)(0)\|^2
+\eps^{\frac{1}{4}}\\
&+\int_0^t\|\xi_{tx}(\tau)\|^2d\tau+\eta_4\int_0^t\|\xi_{txx}(\tau)\|^2d\tau)\\
&\leq C(v_+,\theta_-)(\|\phi_0,\psi_0,\xi_0\|_2^2+\|(\phi_t,\psi_t,\xi_t)(0)\|^2
+\eps^{\frac{1}{4}}\\
&+(\va_0+\eps^{\frac{1}{4}})\int_0^t\|(\phi_{xx},\psi_{xx})(\tau)\|^2
+\eta_4\int_0^t\|\xi_{txx}(\tau)\|^2d\tau),
\end{aligned}
\end{equation}
where $\eta_4$ is a small positive constant to be determined later.
Finally, because of the sufficient small constants $\va_0,\eps$ again,  we see that
\begin{equation}\label{3.74-3}
\begin{aligned}
&\int_0^t\|(\phi_{xx},\psi_{xx})(\tau)\|^2d\tau\\
&\leq C(v_+,\theta_-)(\|\phi_0,\psi_0,\xi_0\|_2^2+\|(\phi_t,\psi_t,\xi_t)(0)\|^2
+\eps^{\frac{1}{4}}+\eta_4\int_{0}^{t}\|\xi_{txx}(\tau)\|^2d\tau).
\end{aligned}
\end{equation}
That is (\ref{3.67}).

Now combining Lemma 3.1-Lemma 3.6, we have following results:

\begin{equation}\label{3.75}
\begin{aligned}
&\|(\phi,\psi,\xi)(t)\|_2^2+\|(\phi_t,\psi_t,\xi_t)(t)\|^2+\int_{0}^{t}\|(\phi_{x},\psi_{x},\xi_x)(\tau)\|_1^2d\tau\\
&+\int_0^t\|\xi_{tx}(\tau)\|^2+\|\xi_{xxx}(\tau)\|^2d\tau\\
&\leq C(v_+,\theta_-)(\|(\phi_0,\psi_0,\xi_0)\|_2^2+\|(\phi_t,\psi_t,\xi_t)(0)\|^2
+\eps^{\frac{1}{4}}+\eta_4\int_{0}^{t}\|\xi_{txx}(\tau)\|^2d\tau).
\end{aligned}
\end{equation}
In order to close our estimate (\ref{3.75}), we have following Lemma.


\begin{Lemma}Under the same assumptions listed in Proposition 3.2, if $\eps, \va_0,\eta_4$ are suitably small, it holds that
	
	\begin{equation}\label{3.76}
	\begin{aligned}
	&\|(\phi_{tx},\psi_{tx},\xi_{tx})(t)\|^2+\int_0^t\|\xi_{txx}(\tau)\|^2d\tau\\
	&\leq C(v_+,\theta_-)(\|(\phi_0,\psi_0,\xi_0)\|_2^2+\|(\phi_t,\psi_t,\xi_t)(0)\|_1^2
	+\eps^{\frac{1}{4}})
	\end{aligned}
	\end{equation}
\end{Lemma}
Proof. $(\ref{3.2})_{1tx}\times\frac{\tilde{p}}{v}\phi_{tx}, (\ref{3.2})_{2t}\times-\psi_{txx}, (\ref{3.2})_{3t}\times-\frac{\xi_{txx}}{\theta}$ imply that 

\begin{equation}\label{3.77}
\begin{aligned}
&\frac{1}{2}(\frac{\tilde{p}}{v}\phi_{tx}^2+\psi_{tx}^2+\frac{R}{\ga-1}\frac{\xi_{tx}^2}{\theta})_t-(\psi_{tt}\psi_{tx}+\frac{R}{\ga-1}\xi_{tt}\frac{\xi_{tx}}{\theta}+\frac{R}{v}\psi_{tx}\xi_{tx})_x
+(\frac{R}{v})_x\psi_{tx}\xi_{tx}\\
&-(\tilde{u}_x(p-\tilde{p}))_t\frac{\xi_{txx}}{\theta}=F_6-k(\frac{\xi_x}{v})_{tx}\frac{\xi_{txx}}{\theta}-k(\frac{\tilde{\theta}_x}{v})_{tx}\frac{\xi_{txx}}{\theta},
\end{aligned}
\end{equation} 
where $F_6$ is

\begin{equation}\label{3.78}
\begin{aligned}
&F_6=\{\frac{1}{2}(\frac{\tilde{p}}{v})_t\phi_{tx}^2+\frac{1}{2}(\frac{R}{\ga-1})(\frac{1}{\theta})_t\xi_{tx}^2-\frac{R}{\ga-1}(\frac{1}{\theta})_x\xi_{tt}\xi_{tx}+p_t\psi_x\frac{\xi_{txx}}{\theta}\}\\
&+\{[((\frac{R}{v})_{x}\xi_t+(\frac{R}{v})_t\xi_x+(\frac{R}{v})_{tx}\xi)
-((\frac{\tilde{p}}{v})_t\phi_x+(\frac{\tilde{p}}{v})_x\phi_t+(\frac{\tilde{p}}{v})_{tx}\phi)]\psi_{txx}\}\\
&=F_{61}+F_{62}.
\end{aligned}
\end{equation}
Obviously,
\begin{equation}
\begin{aligned}\label{3.78-1}
&|F_{61}|\leq C(v_+,\theta_-)[(\va_0+\eps^{\frac{1}{4}})
(\phi_x^2+\psi_x^2+\xi_x^2+\phi_{xx}^2+\psi_{xx}^2+\xi_{xx}^2+\xi_{tx}^2+\xi_{txx}^2)\\
&+|\tilde{v}_x|^2(\phi^2+\xi^2)+|\tilde{\theta}_{txx}|^2+|\tilde{v}_{xx}\tilde{\theta}_x|^2
+|\tilde{v}_x|^6]
\end{aligned}
\end{equation}
Integrating  $(\ref{3.77})$ over $[0,t]\times\mathbb{R}_+$ and using (\ref{3.75}), (\ref{3.78}),(\ref{3.78-1}),
note that $\psi_{tt}(\tau,0)=\xi_{tt}(\tau,0)=0$, we have 

\begin{equation}\label{3.79}
\begin{aligned}
&\|(\phi_{tx},\psi_{tx},\xi_{tx})(t)\|^2+\int_0^t\|\xi_{txx}(\tau)\|^2d\tau\leq C(v_+,\theta_-)(\|(\phi_0,\psi_0,\xi_0)\|_2^2+\|(\phi_t,\psi_t,\xi_t)(0)\|_1^2
+\eps^{\frac{1}{4}}\\
&+(\va_0+\eps^{\frac{1}{4}}+\eta_4)\int_{0}^{t}\|\xi_{txx}(\tau)\|^2d\tau+|\int_{0}^{t}\frac{R}{v}\psi_{tx}\xi_{tx}(\tau,0)d\tau|+\int_{0}^{t}\int_{\mathbb{R}_+}F_{62}dxd\tau.\\
\end{aligned}
\end{equation}
By using the equations $\phi_{tx}=\psi_{xx}$, (\ref{3.34-4}), (\ref{3.61})  and (\ref{3.75}), $F_{62}$ could be estimated by 

\begin{equation}\label{3.80}
\begin{aligned}
&\int_0^t\int_{\mathbb{R}_+}[((\frac{R}{v})_{x}\xi_t+(\frac{R}{v})_t\xi_x+(\frac{R}{v})_{tx}\xi)
-((\frac{\tilde{p}}{v})_t\phi_x+(\frac{\tilde{p}}{v})_x\phi_t+(\frac{\tilde{p}}{v})_{tx}\phi)]\psi_{txx}dxd\tau\\
&=\int_{0}^{t}\int_{\mathbb{R}_+}[((\frac{R}{v})_{x}\xi_t+(\frac{R}{v})_t\xi_x)\psi_{tx}]_x
-((\frac{R}{v})_{x}\xi_t+(\frac{R}{v})_t\xi_x)_x\psi_{tx}\\
&-[((\frac{\tilde{p}}{v})_t\phi_x+(\frac{\tilde{p}}{v})_x\phi_t+\frac{\tilde{p}_{tx}}{v}+\tilde{p}_t(\frac{1}{v})_x+\tilde{p}_x(\frac{1}{v})_t\phi)\psi_{tx}]_x
+(R\xi-\tilde{p}\phi)(\frac{1}{v})_{tx}\psi_{txx}\\
&+((\frac{\tilde{p}}{v})_t\phi_x+(\frac{\tilde{p}}{v})_x\phi_t+\frac{\tilde{p}_{tx}}{v}+\tilde{p}_t(\frac{1}{v})_x+\tilde{p}_x(\frac{1}{v})_t\phi)_x\psi_{tx}dxd\tau\\
&\leq C(v_+,\theta_-)(\|(\phi_0,\psi_0,\xi_0)\|_2^2+\|(\phi_t,\psi_t,\xi_t)(0)\|^2+\eps^{\frac{1}{4}}\\
&+(\eps^{\frac{1}{4}}+\va_0)\int_0^t\phi_x^2(\tau,0)+\psi_x^2(\tau,0)
+\phi_{xx}^2(\tau,0)d\tau+\eta_4\int_{0}^{t}\|\xi_{txx}(\tau)\|^2d\tau\\
&\leq C(v_+,\theta_-)(\|(\phi_0,\psi_0,\xi_0)\|_2^2+\|(\phi_t,\psi_t,\xi_t)(0)\|^2+\eps^{\frac{1}{4}}+\eta_4\int_{0}^{t}\|\xi_{txx}(\tau)\|^2d\tau,
\end{aligned}
\end{equation}
where 
\begin{equation}\label{3.81}
\begin{aligned}
&\int_0^t\int_{\mathbb{R}_+}(R\xi-\tilde{p}\phi)(\frac{1}{v})_{tx}\psi_{txx}dxd\tau\\
&=\int_0^t\int_{\mathbb{R}_+}[(R\xi-\tilde{p}\phi)(2\frac{v_xv_t}{v^3}
-\frac{\tilde{v}_{tx}}{v^2})\psi_{tx}]_x
-[(R\xi-\tilde{p}\phi)(2\frac{v_xv_t}{v^3}-\frac{\tilde{v}_{tx}}{v^2})]_x\psi_{tx}\\
&+((\frac{\tilde{p}\phi-R\xi}{v^2})\frac{\psi_{xx}^2}{2})_t-(\frac{\tilde{p}\phi-R\xi}{v^2})_t\frac{\psi_{xx}^2}{2}dxd\tau.
\end{aligned}
\end{equation}

For the boundary term in the right hand side of $(\ref{3.79}),$  using above results
$(\ref{3.34-4}), (\ref{3.61}), (\ref{3.75})$ again,  we immediately have

\begin{equation}\label{3.82}
\begin{aligned}
&|\int_0^t\frac{R}{v}\psi_{tx}\xi_{tx}(\tau,0)d\tau|\\
&\leq C(v_+,\theta_-)(\|(\phi_0,\psi_0,\xi_0)\|_2^2+\|(\phi_t,\psi_t,\xi_t)(0)\|^2+\eps^{\frac{1}{4}}\\
&+\int_0^t\phi_x^2(\tau,0)+\phi_{xx}^2(\tau,0) +\|\xi_{tx}(\tau)\|^2d\tau+\eta_4\int_0^t\|\xi_{txx}(\tau)\|^2d\tau)\\
&\leq C(v_+,\theta_-)(\|(\phi_0,\psi_0,\xi_0)\|_2^2+\|(\phi_t,\psi_t,\xi_t)(0)\|^2
+\eps^{\frac{1}{4}}+\eta_4\int_0^t\|\xi_{txx}(\tau)\|^2d\tau)).
\end{aligned}
\end{equation}
Substituting $(\ref{3.80})-(\ref{3.82})$ into $(\ref{3.79})$,  $(\ref{3.79})$ turns to 

\begin{equation}\label{3.83}
\begin{aligned}
&\|(\phi_{tx},\psi_{tx},\xi_{tx})(t)\|^2+\int_0^t\|\xi_{txx}(\tau)\|^2d\tau\\
&\leq C(v_+,\theta_-)(\|(\phi_0,\psi_0,\xi_0)\|_2^2+\|(\phi_t,\psi_t,\xi_t)(0)\|_1^2
+\eps^{\frac{1}{4}}
+(\va_0+\eps^{\frac{1}{4}}+\eta_4)\int_0^t\|\xi_{txx}(\tau)\|^2d\tau).
\end{aligned}
\end{equation}
Choosing $\va_0,\eps,\eta_4$ suitably small, we finally get

\begin{equation}\label{3.84}
\begin{aligned}
&\|(\phi_{tx},\psi_{tx},\xi_{tx})(t)\|^2+\int_0^t\|\xi_{txx}(\tau)\|^2d\tau\\
&\leq C(v_+,\theta_-)(\|(\phi_0,\psi_0,\xi_0)\|_2^2+\|(\phi_t,\psi_t,\xi_t)(0)\|_1^2
+\eps^{\frac{1}{4}})
\end{aligned}
\end{equation}
That is (\ref{3.76}). Combining Lemma 3.1-Lemma 3.7, we can get Proposition 3.2.

\end{document}